\documentclass[12pt]{article}

\usepackage{amsfonts}
\usepackage{amssymb}

\title{Vertex operator algebra analogue of
embedding of $B_4$ into $F_4$}
\author{Ozren Per\v{s}e}
\date{}
\begin{document}
\def \Z{\Bbb Z}
\def \C{\Bbb C}
\def \R{\Bbb R}
\def \Q{\Bbb Q}
\def \N{\Bbb N}
\def \tr{{\rm tr}}
\def \span{{\rm span}}
\def \Res{{\rm Res}}
\def \End{{\rm End}}
\def \E{{\rm End}}
\def \Ind {{\rm Ind}}
\def \Irr {{\rm Irr}}
\def \Aut{{\rm Aut}}
\def \Hom{{\rm Hom}}
\def \mod{{\rm mod}}
\def \ann{{\rm Ann}}
\def \<{\langle}
\def \>{\rangle}
\def \t{\tau }
\def \a{\alpha }
\def \e{\epsilon }
\def \l{\lambda }
\def \L{\Lambda }
\def \g{\gamma}
\def \b{\beta }
\def \om{\omega }
\def \o{\omega }
\def \c{\chi}
\def \ch{\chi}
\def \cg{\chi_g}
\def \ag{\alpha_g}
\def \ah{\alpha_h}
\def \ph{\psi_h}
\def \be{\begin{equation}\label}
\def \ee{\end{equation}}
\def \bl{\begin{lem}\label}
\def \el{\end{lem}}
\def \bt{\begin{thm}\label}
\def \et{\end{thm}}
\def \bp{\begin{prop}\label}
\def \ep{\end{prop}}
\def \br{\begin{rem}\label}
\def \er{\end{rem}}
\def \bc{\begin{coro}\label}
\def \ec{\end{coro}}
\def \bd{\begin{de}\label}
\def \ed{\end{de}}
\def \pf{{\bf Proof. }}
\def \voa{{vertex operator algebra}}

\newtheorem{thm}{Theorem}[section]
\newtheorem{prop}[thm]{Proposition}
\newtheorem{coro}[thm]{Corollary}
\newtheorem{conj}[thm]{Conjecture}
\newtheorem{lem}[thm]{Lemma}
\newtheorem{rem}[thm]{Remark}
\newtheorem{de}[thm]{Definition}
\newtheorem{hy}[thm]{Hypothesis}
\makeatletter \@addtoreset{equation}{section}
\def\theequation{\thesection.\arabic{equation}}
\makeatother \makeatletter

\newcommand{\binom}[2]{{{#1}\choose {#2}}}
    \newcommand{\nno}{\nonumber}
    \newcommand{\lbar}{\bigg\vert}
    \newcommand{\p}{\partial}
    \newcommand{\dps}{\displaystyle}
    \newcommand{\bra}{\langle}
    \newcommand{\ket}{\rangle}
 \newcommand{\res}{\mbox{\rm Res}}
\renewcommand{\hom}{\mbox{\rm Hom}}
  \newcommand{\epf}{\hspace{2em}$\Box$}
 \newcommand{\epfv}{\hspace{1em}$\Box$\vspace{1em}}
\newcommand{\nord}{\mbox{\scriptsize ${\circ\atop\circ}$}}
\newcommand{\wt}{\mbox{\rm wt}\ }

\maketitle
\begin{abstract}
Let $L_{B}(-\frac{5}{2},0)$ (resp. $L_{F}(-\frac{5}{2},0)$) be the
simple vertex operator algebra associated to affine Lie algebra of
type $B_{4}^{(1)}$ (resp. $F_{4}^{(1)}$) with the lowest admissible
half-integer level $-\frac{5}{2}$. We show that
$L_{B}(-\frac{5}{2},0)$ is a vertex subalgebra of
$L_{F}(-\frac{5}{2},0)$ with the same conformal vector, and that
$L_{F}(-\frac{5}{2},0)$ is isomorphic to the extension of
$L_{B}(-\frac{5}{2},0)$ by its only irreducible module other than
itself. We also study the representation theory of
$L_{F}(-\frac{5}{2},0)$, and determine the decompositions of
irreducible weak $L_{F}(-\frac{5}{2},0)$-modules from the category
$\mathcal{O}$ into direct sums of irreducible weak
$L_{B}(-\frac{5}{2},0)$-modules.
\end{abstract}


\footnotetext[2]{ Partially supported by the
Ministry of Science, Education and Sports
of the Republic of Croatia, grant 0037125.}

\section{Introduction}

The embedding of $B_4$ into $F_4$ has been studied on the level of
Lie groups and Lie algebras both in mathematics and physics (cf.
\cite{Ba,GKRS,PR}). It is particularly interesting that simple Lie
algebra of type $F_4$ contains three copies of simple Lie algebra of
type $B_4$ as Lie subalgebras.

Let $\frak g _{F}$ be the simple Lie algebra of type $F_{4}$, and
$\frak g _{B}$ its Lie subalgebra of type $B_{4}$. The decomposition
of $\frak g _{F}$ into a direct sum of irreducible $\frak g
_{B}$-modules is:
\begin{eqnarray} \label{rel.intro}
\frak g _{F} \cong \frak g _{B} \oplus V_{B}(\bar{\omega} _{4}),
\end{eqnarray}
where $V_{B}(\bar{\omega} _{4})$ is the irreducible highest weight
$\frak g _{B}$-module whose highest weight is the fundamental weight
$\bar{\omega} _{4}$ for $\frak g _{B}$.

We want to study the analogue of relation (\ref{rel.intro}) in the
case of vertex operator algebras associated to affine Lie algebras.
Let $\hat{\frak g}_{B}$ be the affine Lie algebra of type
$B_{4}^{(1)}$, and $\hat{\frak g}_{F}$ the affine Lie algebra of
type $F_{4}^{(1)}$. For $k \in \C$, denote by $L_{B}(k,0)$ and
$L_{F}(k,0)$ simple vertex operator algebras associated to
$\hat{\frak g}_{B}$ and $\hat{\frak g}_{F}$ of level $k$,
respectively, with conformal vectors obtained from Segal-Sugawara
construction.

We want to study a certain level $k$, such that $L_{B}(k,0)$ is a
vertex subalgebra of $L_{F}(k,0)$ with the same conformal vector.
The equality of conformal vectors implies the equality of
corresponding central charges
\begin{eqnarray*}
\frac{k \dim \frak g _{B}}{k+h_{B}^{\vee}}= \frac{k \dim \frak g
_{F}}{k+h_{F}^{\vee}}.
\end{eqnarray*}
This is a linear equation in $k$, and its only solution is
$k=-\frac{5}{2}$. This suggests the study of vertex operator
algebras $L_{B}(-\frac{5}{2},0)$ and $L_{F}(-\frac{5}{2},0)$.

Vertex operator algebra $L_{B}(-\frac{5}{2},0)$ is a special case of
vertex operator algebras associated to type $B$ affine Lie algebras
of admissible half-integer levels, studied in \cite{P}. Admissible
weights are a class of weights defined in \cite{KW1} and \cite{KW2},
that includes dominant integral weights. The character formula for
admissible modules is given in these papers, which inspires the
study of vertex operator algebras with admissible levels. In case of
non-integer admissible levels, these vertex operator algebras are
not rational in the sense of Zhu \cite{Z}, but in all known cases,
they have properties similar to rationality in the category of weak
modules that are in category $\mathcal{O}$ as modules for
corresponding affine Lie algebra (cf. \cite{A1,AM,DLM1}). Admissible
modules for affine Lie algebras were also recently studied in
\cite{A3,FM,GPW,W}.

In case of affine Lie algebra $\hat{\frak g}_{B}$ of type
$B_{4}^{(1)}$, levels $n-\frac{7}{2}$ are studied in \cite{P}, for
$n \in \N$. For certain subcategories of the category of weak
$L_{B}(n-\frac{7}{2},0)$-modules, the irreducible objects are
classified and semisimplicity is proved. It follows that in case
$n=1$, only irreducible $L_{B}(-\frac{5}{2},0)$-modules are
$L_{B}(-\frac{5}{2},0)$ and $L_{B}(-\frac{5}{2}, \bar{\omega}
_{4})$.

Another motivation for studying this case is a problem of extension
of a vertex operator algebra by a module. This problem was studied
in \cite{DLM2}, where all simple current extensions of vertex
operator algebras associated to affine Lie algebras (except for
$E_{8}^{(1)}$) of positive integer levels $k$ are constructed. These
extensions have a structure of abelian intertwining algebra, defined
and studied in \cite{DL}. In some special cases, these extensions
have a structure of a vertex operator algebra. In our case of
admissible non-integer level $k=-\frac{5}{2}$, it is natural to
consider the extension of vertex operator algebra
$L_{B}(-\frac{5}{2},0)$ by its only irreducible module other then
itself, $L_{B}(-\frac{5}{2}, \bar{\omega} _{4})$.

In this paper we study vertex operator algebras
$L_{F}(n-\frac{7}{2},0)$, for $n \in \N$. We show that levels
$n-\frac{7}{2}$ are admissible for $\hat{\frak g}_{F}$. Results on
admissible modules from \cite{KW1} imply that
$L_{F}(n-\frac{7}{2},0)$ is a quotient of the generalized Verma
module by the maximal ideal generated by one singular vector. We
determine the formula for that singular vector and using that, we
show that $L_{F}(n-\frac{7}{2},0)$ contains three copies of
$L_{B}(n-\frac{7}{2},0)$ as vertex subalgebras.

In the case $n=1$, we show that these three copies of
$L_{B}(-\frac{5}{2},0)$ have the same conformal vector as
$L_{F}(-\frac{5}{2},0)$, which along with results from \cite{P}
implies that
\begin{eqnarray*}
L_{F}(-\frac{5}{2},0) \cong L_{B}(-\frac{5}{2},0) \oplus
L_{B}(-\frac{5}{2}, \bar{\omega} _{4}),
\end{eqnarray*}
which is a vertex operator algebra analogue of relation
(\ref{rel.intro}). This also implies that the extension of vertex
operator algebra $L_{B}(-\frac{5}{2},0)$ by its module
$L_{B}(-\frac{5}{2}, \bar{\omega} _{4})$ has a structure of vertex
operator algebra, isomorphic to $L_{F}(-\frac{5}{2},0)$.

Furthermore, we study the category of weak
$L_{F}(-\frac{5}{2},0)$-modules that are in category $\mathcal{O}$
as $\hat{\frak g}_{F}$-modules. Using three copies of
$L_{B}(-\frac{5}{2},0)$ contained in $L_{F}(-\frac{5}{2},0)$,
methods from \cite{MP,A2,AM} and results from \cite{P}, we obtain
the classification of irreducible objects in that category. Using
results from \cite{KW2} we prove semisimplicity of that category. We
also determine decompositions of irreducible weak
$L_{F}(-\frac{5}{2},0)$-modules from category $\mathcal{O}$ into
direct sums of irreducible weak $L_{B}(-\frac{5}{2},0)$-modules. It
turns out that these direct sums are finite.

We also study the category of $L_{F}(n-\frac{7}{2},0)$-modules with
higher admissible half-integer levels $n-\frac{7}{2}$, for $n \in
\N$. Using results from \cite{P}, we classify irreducible objects in
that category and prove semisimplicity of that category.

The author expresses his gratitude to Professor Dra\v{z}en Adamovi\'{c} for
his helpful advice and valuable suggestions.

\section{Vertex operator algebras associated to affine Lie algebras}
\label{sect.prelim}

This section is preliminary. We recall some necessary definitions
and fix the notation. We review certain results about vertex
operator algebras and corresponding modules. The emphasis is on
the class of vertex operator algebras associated to affine Lie
algebras, because we study special cases in that class.

\subsection{Vertex operator algebras and modules}

Let $(V, Y, {\bf 1}, \omega)$ be a vertex operator algebra (cf.
\cite{Bo}, \cite{FHL}, \cite{FLM} and \cite{LL}). Specially, the
triple  $(V, Y, {\bf 1})$ carries the structure of a vertex algebra.

A {\it vertex subalgebra} of vertex algebra $V$ is a subspace
$U$ of $V$ such that ${\bf 1} \in U$ and $ Y(a,z)U \subseteq U[[z,z^{-1}]] $
for any $a \in U$.
Suppose that $(V, Y, {\bf 1}, \omega)$ is a vertex operator algebra
and $(U, Y, {\bf 1}, \omega ')$ a vertex subalgebra of $V$,
that has a structure of vertex operator algebra. We say that $U$ is
a {\it vertex operator subalgebra} of $V$ if $\omega ' = \omega$.

An {\it ideal} in a vertex operator algebra $V$ is a subspace $I$
of $V$ satisfying $ Y(a,z)I \subseteq I[[z,z^{-1}]] $
for any $a \in V$. Given an ideal $I$ in $V$,
such that ${\bf 1} \notin I$, $\omega \notin I$, the
quotient $V/I$ admits a natural vertex operator algebra structure.

Let $(M, Y_{M})$ be a weak module for a vertex operator algebra
$V$ (cf. \cite{L}).
A {\it ${\Z}_{+}$-graded weak} $V$-module (\cite{FZ}) is a weak $V$-module
$M$ together with a ${\Z}_{+}$-gradation
$M=\oplus_{n=0}^{\infty}M(n)$ such that
\[
a_{m}M(n)\subseteq M(n+r-m-1)\;\;\;\mbox{for }a \in
V_{(r)},m,n,r\in {\Z},
\]
where $M(n)=0$ for $n < 0$ by definition.

A weak $V$-module $M$ is called a {\it $V$-module} if $L(0)$
acts semisimply on $M$ with the decomposition into
$L(0)$-eigenspaces $M=\oplus_{\alpha\in {\C}}M_{(\alpha)}$ such
that for any $\alpha \in {\C}$, $\dim M_{(\alpha)}<\infty$
and $M_{(\alpha+n)}=0$ for $n\in {\Z}$ sufficiently small.

\subsection{Zhu's $A(V)$ theory}

Let $V$ be a vertex operator algebra. Following \cite{Z}, we
define bilinear maps $* :  V \times V
\to V$ and $\circ  :  V \times V \to V$ as follows.
For any homogeneous $a \in V$ and for any $b \in V$, let
\begin{eqnarray*}
a \circ b=\Res_{z}\frac{(1+z)^{{\rm wt} a}}{z^{2}}Y(a,z)b  \\
a * b=\Res_{z}\frac{(1+z)^{{\rm wt} a}}{z}Y(a,z)b
\end{eqnarray*}
and extend to $V \times V \to V$ by linearity. Denote by
$O(V)$ the linear span of elements of the form $a \circ b$, and by
$A(V)$ the quotient space $V/O(V)$. For $a \in V$, denote by
$[a]$ the image of $a$ under the projection of $V$ onto
$A(V)$. The multiplication $*$ induces the multiplication
on $A(V)$ and $A(V)$ has a structure of an associative algebra.

\begin{prop}[{\cite[Proposition 1.4.2]{FZ}}] Let $I$ be an ideal of $V$. Assume
${\bf 1} \notin I$, $\omega \notin I$. Then the associative
algebra $A(V/I)$ is isomorphic to $A(V)/A(I)$, where
$A(I)$ is the image of $I$ in $A(V)$.
\end{prop}

For any homogeneous $a \in V$ we define $o(a)=a_{{\rm wt} a-1}$
and extend this map linearly to $V$.

\begin{prop} [{\cite[Theorem 2.1.2, Theorem 2.2.1]{Z}}]
\item[(a)] Let $M=\oplus_{n=0}^{\infty}M(n)$ be a ${\Z}_{+}$-graded
weak $V$-module. Then $M(0)$ is an $A(V)$-module defined as
follows:
\[
[a].v=o(a)v,
\]
for any $a \in V$ and $v \in M(0)$.
\item[(b)] Let $U$ be an $A(V)$-module. Then there exists a ${\Z}_{+}$-graded
weak $V$-module $M$ such that the $A(V)$-modules $M(0)$ and $U$ are
isomorphic.
\end{prop}

\begin{prop}[{\cite[Theorem 2.2.2]{Z}}] \label{t.1.3.5}
The equivalence classes of the irreducible $A(V)$-modules and
the equivalence classes of the irreducible ${\Z}_{+}$-graded weak $V$-modules
are in one-to-one correspondence.
\end{prop}

\subsection{Modules for affine Lie algebras}

Let ${\frak g}$ be a simple Lie algebra over ${\C}$ with a
triangular decomposition \linebreak ${\frak g}={\frak n}_{-} \oplus
{\frak h} \oplus {\frak n}_{+}$. Let $\Delta$ be the root system of
$({\frak g}, {\frak h})$, $\Delta_{+}\subset \Delta$
the set of positive roots, $\theta$ the highest root and
$(\cdot, \cdot): {\frak g}\times {\frak g}\to {\C}$ the Killing form,
normalized by the condition $(\theta, \theta)=2$.

The affine Lie algebra $\hat{\frak g}$ associated to ${\frak g}$
is the vector space ${\frak g}\otimes {\C}[t, t^{-1}] \oplus {\C}c
$ equipped with the usual bracket operation and canonical central
element~$c$ (cf. \cite{K}).
Let $h^{\vee}$ be the dual Coxeter number of $\hat{\frak g}$.
Let $\hat{\frak g}=\hat{\frak n}_{-} \oplus \hat{\frak h}
\oplus \hat{\frak n}_{+}$ be the corresponding
triangular decomposition of $\hat{\frak g}$.
Denote by $\hat{\Delta}$ the set of roots of $\hat{\frak g}$,
by $\hat{\Delta}^{\mbox{\scriptsize{re}}}$
the set of real roots of $\hat{\frak g}$, and by
$\alpha ^{\vee}$ denote the coroot of
a real root $\alpha \in \hat{\Delta}^{\mbox{\scriptsize{re}}}$.

For every weight $\lambda\in \hat{{\frak h}}^{*}$, denote by
$M( \lambda)$ the Verma module for $\hat{\frak g}$ with highest
weight $\lambda$, and by $L( \lambda)$ the irreducible
$\hat{\frak g}$-module with highest
weight~$\lambda$.

Let $U$ be a ${\frak g}$-module, and let $k\in {\C}$. Let
$\hat{\frak g}_{+}={\frak g}\otimes t{\C}[t]$ act trivially on $U$
and $c$ as scalar $k$. Considering $U$ as a ${\frak g}\oplus {\C}c
\oplus \hat{\frak g}_{+}$-module, we have the induced $\hat{\frak
g}$-module (so called {\it generalized Verma module})
\[
N(k,U)=U(\hat{\frak g})\otimes_{U({\frak g}\oplus {\C}c
\oplus \hat{\frak g}_{+})} U.
\]

For a fixed $\mu \in {\frak h}^{*}$, denote by $V(\mu)$ the
irreducible highest weight \linebreak ${\frak g}$-module with highest weight
$\mu$.

We shall use the notation $N(k, \mu)$ to denote the $\hat{\frak
g}$-module $N(k,V(\mu))$. Denote by $J(k, \mu)$ the
maximal proper submodule of $N(k, \mu)$ and by $L(k, \mu)$
the corresponding irreducible quotient $N(k, \mu)/J(k, \mu)$.

\subsection{Admissible modules for affine Lie algebras}

Let $\hat{\Delta}^{\vee \mbox{\scriptsize{re}}}$\
(resp. $\hat{\Delta}^{\vee \mbox{\scriptsize{re}}}_{+}$)
$\subset \hat{\frak h}$ be the set of real
(resp. positive real) coroots of $ \hat{\frak g}$. Fix $\lambda \in \hat{\frak h}^*$.
Let $\hat{\Delta}^{\vee \mbox{\scriptsize{re}}}_{\lambda}=
\{\alpha\in \hat{\Delta}^{\vee \mbox{\scriptsize{re}}} \ |\ \langle\lambda,\alpha\rangle\in{\Z} \}$,
$\hat{\Delta}^{\vee \mbox{\scriptsize{re}}}_{\lambda +}=\hat{\Delta}^{\vee \mbox{\scriptsize{re}}}_{\lambda}
\cap \hat{\Delta}^{\vee \mbox{\scriptsize{re}}}_{+}$,
$\hat{\Pi}^{\vee}$ the set of simple coroots in
$\hat{\Delta}^{\vee \mbox{\scriptsize{re}}}$ and
$\hat{\Pi}^{\vee}_{\lambda }=
\{\alpha \in \hat{\Delta}^{\vee \mbox{\scriptsize{re}}}_{\lambda +}\ \vert \ \alpha$
not equal to a sum of several coroots from $\hat{\Delta}^{\vee \mbox{\scriptsize{re}}}_{\lambda +} \}$.
Define $\rho$ in the usual way, and denote by $w. \lambda$
the "shifted" action of an element $w$ of the Weyl group of
$ \hat{\frak g}$.

Recall that a weight  $\lambda \in \hat{\frak h} ^* $
is called {\it admissible} (cf. \cite{KW1}, \cite{KW2} and \cite{W})
if the following properties are satisfied:
\begin{eqnarray*}
& &\langle\lambda + \rho,\alpha\rangle \notin -{\Z}_+ \mbox{ for all }
\alpha \in \hat{\Delta}^{\vee \mbox{\scriptsize{re}}}_{+}, \\
& &{\Q} \hat{\Delta}^{\vee \mbox{\scriptsize{re}}}_{\lambda}={\Q} \hat{\Pi}^{\vee} .
\end{eqnarray*}
The irreducible $\hat{\frak g}$-module $L(\lambda)$
is called admissible if the weight $\lambda \in \hat{\frak h} ^* $
is admissible.

We shall use the following results of V. Kac and M. Wakimoto:
\begin{prop}[{\cite[Corollary 2.1]{KW1}}] \label{t.KW1}
Let $\lambda $ be an admissible weight. Then
\[
L(\lambda) \cong \frac {M(\lambda)}{\sum_{\alpha \in \hat{\Pi}^{\vee}_{\lambda }}U (
\hat{\frak g}) v^{\alpha}}\ ,
\]
where $v^{\alpha}\in M(\lambda)$ is a singular vector of weight $r_{\alpha}.
\lambda$, the highest weight vector of $M(r_{\alpha}.\lambda)=\ U(\hat{\frak g})
v^{\alpha}\subset M(\lambda)$.
\end{prop}

\begin{prop}[{\cite[Theorem 4.1]{KW2}}] \label{t.KW2}
Let $M$ be a $\hat{\frak g}$-module from the category $\mathcal{O}$ such
that for any irreducible subquotient $L(\nu )$ the weight $\nu $ is
admissible. Then $M$ is completely reducible.
\end{prop}

\subsection{Vertex operator algebras $N(k,0)$ and $L(k,0)$, for
\\ $k \neq - h^{\vee}$}

Since $V(0)$ is the one-dimensional
trivial ${\frak g}$-module, it can be identified with~${\C}$.
Denote by ${\bf 1}=1 \otimes 1 \in N(k,0)$.
We note that $N(k,0)$ is spanned by the elements of the form
$x_{1}(-n_{1}-1)\cdots x_{m}(-n_{m}-1){\bf 1}$, where  $x_{1}, \dots, x_{m}\in
{\frak g}$ and $n_{1}, \dots, n_{m}\in {\Z}_{+}$, with $x(n)$ denoting
the representation image of
$x\otimes t^{n}$ for $x\in{\frak g}$ and $n\in {\Z}$.
Vertex operator map
$Y(\cdot, z): N(k,0) \to  (\mbox{\rm End}\;N(k,0))[[z, z^{-1}]]$
is uniquely determined by defining $Y({\bf 1}, z)$ to be the
identity operator on $N(k,0)$ and
\[
Y(x(-1){\bf 1}, z)=\sum_{n\in {\Z}}x(n)z^{-n-1},
\]
for $x\in {\frak g}$.
In the case that $k\ne - h^{\vee}$, $N(k,0)$ has a conformal
vector
\begin{eqnarray} \label{rel.Virasoro}
\omega=\frac{1}{2(k+h^{\vee})}\sum_{i=1}^{\dim {\frak g}}
a^{i}(-1)b^{i}(-1){\bf 1},
\end{eqnarray}
where $\{a^{i}\}_{i=1, \dots, \dim {\frak g}}$ is an arbitrary
basis of ${\frak g}$, and
$\{b^{i}\}_{i=1, \dots, \dim {\frak g}}$  the corresponding dual
basis of ${\frak g}$ with respect to the form $(\cdot, \cdot)$.
We have the following result from \cite{FZ}, \cite{L}:

\begin{prop}
If $k\ne -h^{\vee}$, the quadruple
$(N(k,0), Y, {\bf 1}, \omega)$ defined above is a vertex operator algebra.
\end{prop}

The associative algebra $A(N(k,0))$ is identified in
next proposition:

\begin{prop} [{\cite[Theorem 3.1.1]{FZ}}]
The associative algebra $A(N(k,0))$
is canonically isomorphic to $U (\frak g ) $.
The isomorphism is given by $F:A(N(k,0)) \to U (\frak g )$
\[
F([x_1(-n_1 -1)\cdots x_m(-n_m -1){\bf 1}])= (-1)^{n_1+\cdots +n_m}
x_m \cdots x_1,
\]
for any $x_1, \ldots ,x_m \in \frak g$ and any $n_1, \ldots ,n_m \in {\Z}_{+}$.
\end{prop}

For any $\mu \in {\frak h}^{*}$, $N(k,\mu)$ is a $\Z _{+}$-graded
weak $N(k,0)$-module, and $L(k,\mu)$
is an irreducible $\Z _{+}$-graded
weak $N(k,0)$-module. Denote by $v_{k, \mu}$ the highest weight
vector of $L(k,\mu)$. Then the lowest conformal weight of
$L(k,\mu)$ is given by relation
\begin{eqnarray} \label{rel.lowest.conf.w}
L(0)v_{k, \mu}=\frac{(\mu, \mu + 2 \bar{\rho})}{2(k+h^{\vee})}v_{k,
\mu},
\end{eqnarray}
where $\bar{\rho}$ is the sum of fundamental weights of $\frak g$.

Since every $\hat{\frak g}$-submodule of $N(k,0)$ is also an ideal
in the vertex operator algebra $N(k,0)$, it follows that
$L(k,0)$ is a vertex operator algebra, for any $k\ne -h^{\vee}$.
The associative algebra $A(L(k,0))$ is identified in
next proposition, in the case when maximal $\hat{\frak g}$-submodule of
$N(k,0)$ is generated by one singular vector.

\begin{prop} \label{p.1.5.5}
Assume that the maximal $\hat{\frak g}$-submodule of
$N(k,0)$ is generated by a singular vector, i.e.
$J(k,0)=U(\hat{\frak g})v.$
Then
\[
A(L(k,0)) \cong \frac{U(\frak g)}{I},
\]
where $I$ is the two-sided ideal of $U(\frak g)$
generated by $u=F([v])$. \\ \noindent
Let $U$ be a $\frak g$-module. Then $U$ is an
$A(L(k,0))$-module if and only if $IU=0$.
\end{prop}

\subsection{Modules for associative algebra $A(L(k,0))$}
\label{subsec.3.3}

In this subsection we present the method from \cite{MP},
\cite{A2}, \cite{AM} for
classification of irreducible $A(L(k,0))$-modules
from the category $\mathcal{O}$ by solving certain systems of
polynomial equations. We assume that the maximal
$\hat{\frak g}$-submodule of
$N(k,0)$ is generated by a singular vector $v$.

Denote by $_L$ the adjoint action of  $U(\frak g)$ on
$U(\frak g)$ defined by $ X_Lf=[X,f]$ for $X \in
\frak g$ and $f \in U(\frak g)$. Let $R$ be a $U(\frak g)$-submodule
of $U(\frak g)$ generated by the vector $u=F([v])$ under the adjoint
action. Clearly, $R$ is an irreducible highest weight $U(\frak
g)$-module with the highest weight vector $u$. Let $R_{0}$ be the zero-weight
subspace of $R$. The next proposition follows from
\cite[Proposition 2.4.1]{A2}, \cite[Lemma 3.4.3]{AM}:

\begin{prop}
Let $V(\mu)$ be an irreducible highest weight $U(\frak g)$-module
with the highest weight vector $v_{\mu}$, for $\mu \in {\frak h}^{*}$.
The following statements are equivalent:
\item[(1)] $V(\mu)$ is an $A(L(k,0))$-module,
\item[(2)] $RV(\mu)=0$,
\item[(3)] $R_{0}v_{\mu}=0.$
\end{prop}

Let $r \in R_{0}$. Clearly there exists the unique polynomial
$p_{r} \in S( \frak h)$ such that
\[
rv_{\mu}=p_{r}(\mu)v_{\mu}.
\]
Set $ {\mathcal P}_{0}=\{ \ p_{r} \ \vert \ r \in R_{0} \}.$
We have:

\begin{coro} \label{c.1.7.2} There is one-to-one correspondence between
\item[(1)] irreducible $A(L(k,0))$-modules
from the category $\mathcal{O}$,
\item[(2)] weights $\mu \in {\frak h}^{*}$ such that $p(\mu)=0$
for all $p \in {\mathcal P}_{0}$.
\end{coro}

\section{Simple Lie algebras of type $B_{4}$ and $F_{4}$}
\label{sect.simple}

Let
\[
\Delta _{F} =\{ \pm \epsilon_{i} \, \vert \ 1 \leq i \leq 4 \} \cup
\{ \pm \epsilon_{i} \pm \epsilon_{j} \, \vert \ 1 \leq i < j \leq 4 \}
\cup \{ \frac{1}{2}( \pm \epsilon_{1} \pm \epsilon_{2}
\pm \epsilon_{3} \pm \epsilon_{4}) \}
\]
be the root system of type $F_{4}$. Fix the set of positive roots
$\Delta_{F}^{+}=\{ \epsilon_{i}, \, \vert \ 1 \leq i \leq 4 \} \cup
\{ \epsilon_{i} \pm \epsilon_{j} \, \vert \ 1 \leq i < j \leq 4 \} \cup
\{ \frac{1}{2}( \epsilon_{1} \pm \epsilon_{2}
\pm \epsilon_{3} \pm \epsilon_{4}) \}$.
Then the simple roots are
\[
\alpha_{1}= \epsilon_{2} - \epsilon_{3},
\alpha_{2}= \epsilon_{3} - \epsilon_{4},
\alpha_{3}= \epsilon_{4},
\alpha_{4}= \frac{1}{2}( \epsilon_{1} - \epsilon_{2}
- \epsilon_{3} - \epsilon_{4}).
\]
The highest root is $\theta = \epsilon_{1} + \epsilon_{2}$.

The subset $\Delta _{B} \subseteq \Delta _{F}$, defined by
\[
\Delta _{B} =\{ \pm \epsilon_{i} \, \vert \ \ 1 \leq i \leq 4 \} \cup
\{ \pm \epsilon_{i} \pm \epsilon_{j} \, \vert \ 1 \leq i < j \leq 4 \}
\]
is a root system of type $B_{4}$. Clearly, $\theta \in \Delta _{B}$.
The corresponding simple roots are
\[
\beta_{1}= \epsilon_{1} - \epsilon_{2},
\beta_{2}= \epsilon_{2} - \epsilon_{3},
\beta_{3}= \epsilon_{3} - \epsilon_{4},
\beta_{4}= \epsilon_{4}.
\]
Furthermore, subsets $\Delta _{B}', \Delta _{B}''  \subseteq \Delta _{F}$
defined by
\begin{eqnarray*}
&& \Delta _{B}' =\{ \pm \frac{1}{2}( \epsilon_{1} - \epsilon_{2}
- \epsilon_{3} - \epsilon_{4}), \pm \frac{1}{2}( \epsilon_{1} + \epsilon_{2}
- \epsilon_{3} + \epsilon_{4}), \pm \frac{1}{2}( \epsilon_{1} - \epsilon_{2}
+ \epsilon_{3} + \epsilon_{4}),  \\
&& \qquad \pm \frac{1}{2}( \epsilon_{1} + \epsilon_{2}
+ \epsilon_{3} - \epsilon_{4})  \} \cup
\{ \pm \epsilon_{i} \pm \epsilon_{j}
\, \vert \ 1 \leq i < j \leq 4 \}, \\
&& \Delta _{B}'' =\{ \pm \frac{1}{2}( \epsilon_{1} - \epsilon_{2}
- \epsilon_{3} + \epsilon_{4}), \pm \frac{1}{2}( \epsilon_{1} - \epsilon_{2}
+ \epsilon_{3} - \epsilon_{4}), \pm \frac{1}{2}( \epsilon_{1} - \epsilon_{2}
- \epsilon_{3} + \epsilon_{4}),  \\
&& \qquad \pm \frac{1}{2}( \epsilon_{1} + \epsilon_{2}
+ \epsilon_{3} + \epsilon_{4})  \} \cup
\{ \pm \epsilon_{i} \pm \epsilon_{j}
\, \vert \ 1 \leq i < j \leq 4 \}
\end{eqnarray*}
are also root systems of type $B_{4}$. The corresponding simple
roots are
\[
\beta_{1}'= \epsilon_{3} - \epsilon_{4},
\beta_{2}'= \epsilon_{2} - \epsilon_{3},
\beta_{3}'= \epsilon_{3} + \epsilon_{4},
\beta_{4}'= \frac{1}{2}( \epsilon_{1} - \epsilon_{2}
- \epsilon_{3} - \epsilon_{4})
\]
for $\Delta _{B}'$ and
\[
\beta_{1}''= \epsilon_{3} + \epsilon_{4},
\beta_{2}''= \epsilon_{2} - \epsilon_{3},
\beta_{3}''= \epsilon_{3} - \epsilon_{4},
\beta_{4}''= \frac{1}{2}( \epsilon_{1} - \epsilon_{2}
- \epsilon_{3} + \epsilon_{4})
\]
for $\Delta _{B}''$.

Let $\frak g _{F}$ be the simple Lie algebra associated to
the root system $\Delta _{F}$, and $\frak g _{B}$, $\frak g _{B}'$ and
$\frak g _{B}''$ its Lie subalgebras
associated to root systems $\Delta _{B}$, $\Delta _{B}'$ and $\Delta _{B}''$,
respectively. The isomorphism $\pi '$ of root systems $\Delta _{B}$ and $\Delta _{B}'$
defined by
\[
\pi '(\beta_{i})= \beta_{i}', \quad i=1,2,3,4
\]
induces an isomorphism of Lie algebras $\pi ' : \frak g _{B} \to
\frak g _{B}'$. Similarly, the isomorphism $\pi ''$ of root systems
$\Delta _{B}$ and $\Delta _{B}''$ defined by
\[
\pi ''(\beta_{i})= \beta_{i}'', \quad i=1,2,3,4
\]
induces an isomorphism of Lie algebras $\pi '' : \frak g _{B} \to
\frak g _{B}''$. Thus, we have three ways of embedding simple Lie
algebra of type $B_{4}$ into simple Lie algebra of type $F_{4}$.

Let $e_{i},f_{i},h_{i}$, $1 \leq i \leq 4$ be the Chevalley generators of $\frak g _{F}$.
Fix the root vectors for positive roots of $\frak g _{F}$:
\begin{eqnarray} \label{rel.kor.vekt.}
&& e_{\frac{1}{2}( \epsilon_{1} - \epsilon_{2}
- \epsilon_{3} + \epsilon_{4})}=[e_{3},e_{4}], \
e_{\epsilon_2 - \epsilon_4}=[e_{2},e_{1}], \
e_{\epsilon_3}=[e_{2},e_{3}], \nonumber \\
&& e_{\frac{1}{2}( \epsilon_{1} - \epsilon_{2}
+ \epsilon_{3} - \epsilon_{4})}=[e_{\epsilon_3}
,e_{\frac{1}{2}( \epsilon_{1} - \epsilon_{2}
- \epsilon_{3} + \epsilon_{4})}], \
e_{\epsilon_2}=[e_{\epsilon_2 - \epsilon_4},e_{3}], \
e_{\epsilon_3 + \epsilon_4}=\frac{1}{2}[e_{\epsilon_3},e_{3}], \nonumber \\
&& e_{\frac{1}{2}( \epsilon_{1} + \epsilon_{2}
- \epsilon_{3} - \epsilon_{4})}=[e_{\frac{1}{2}( \epsilon_{1} - \epsilon_{2}
+ \epsilon_{3} - \epsilon_{4})}, e_{1}], \
e_{\frac{1}{2}( \epsilon_{1} - \epsilon_{2}
+ \epsilon_{3} + \epsilon_{4})}=[e_{\frac{1}{2}( \epsilon_{1} - \epsilon_{2}
- \epsilon_{3} + \epsilon_{4})}, e_{3}], \nonumber \\
&& e_{\epsilon_2 + \epsilon_4}=\frac{1}{2}[e_{\epsilon_2},e_{3}],
\ e_{\frac{1}{2}( \epsilon_{1} + \epsilon_{2}
- \epsilon_{3} + \epsilon_{4})}=[e_{\frac{1}{2}( \epsilon_{1} + \epsilon_{2}
- \epsilon_{3} - \epsilon_{4})}, e_{3}], \nonumber \\
&& e_{\epsilon_1 - \epsilon_2}=\frac{1}{2}[e_{\frac{1}{2}( \epsilon_{1} - \epsilon_{2}
+ \epsilon_{3} - \epsilon_{4})}, e_{\frac{1}{2}( \epsilon_{1} - \epsilon_{2}
- \epsilon_{3} + \epsilon_{4})}], \
e_{\epsilon_2 +
\epsilon_3}=\frac{1}{2}[e_{\epsilon_3},e_{\epsilon_2}], \nonumber \\
&& e_{\frac{1}{2}( \epsilon_{1} + \epsilon_{2}
+ \epsilon_{3} - \epsilon_{4})}=[e_{\frac{1}{2}( \epsilon_{1} + \epsilon_{2}
- \epsilon_{3} + \epsilon_{4})}, e_{2}], \
e_{\epsilon_1 - \epsilon_3}=[e_{\epsilon_1 - \epsilon_2},e_{1}], \nonumber \\
&& e_{\frac{1}{2}( \epsilon_{1} + \epsilon_{2}
+ \epsilon_{3} + \epsilon_{4})}=[e_{3}, e_{\frac{1}{2}( \epsilon_{1} + \epsilon_{2}
+ \epsilon_{3} - \epsilon_{4})}],\
e_{\epsilon_1 - \epsilon_4}=[e_{\epsilon_1 - \epsilon_3},e_{2}],\
e_{\epsilon_1}=[e_{\epsilon_1 - \epsilon_4},e_{3}], \nonumber \\
&& e_{\epsilon_1 +
\epsilon_4}=\frac{1}{2}[e_{\epsilon_1},e_{3}], \
e_{\epsilon_1 +
\epsilon_3}=\frac{1}{2}[e_{\epsilon_1},e_{\epsilon_3}], \
e_{\epsilon_1 +
\epsilon_2}=\frac{1}{2}[e_{\epsilon_2},e_{\epsilon_1}].
\end{eqnarray}
To define root vectors for negative roots of $\frak g _{F}$,
put $f_{\alpha + \beta}=-c_{\alpha , \beta}[f_{\alpha},f_{\beta}]$
if $e_{\alpha + \beta}=c_{\alpha , \beta}[e_{\alpha},e_{\beta}]$.
Denote by $h_{\alpha}= \alpha ^{\vee}= [e_{\alpha},f_{\alpha}]$
coroots, for any positive root $\alpha \in \Delta _{F}^{+}$.
Let ${\frak g}_{F}={\frak n}_{-} \oplus
{\frak h} \oplus {\frak n}_{+}$ be the corresponding
triangular decomposition of ${\frak g}_{F}$.

One can easily check that $\pi '(e_{\alpha})=e_{\pi '(\alpha)}$,
$\pi '(f_{\alpha})=f_{\pi '(\alpha)}$,
and $\pi ''(e_{\alpha})=e_{\pi ''(\alpha)}$,
$\pi ''(f_{\alpha})=f_{\pi ''(\alpha)}$ for all $\alpha \in \Delta
_{B}^+$.

Denote by $P_{+}^{F}$ the set of dominant integral weights for
$\frak g_{F}$ and by  $P_{+}^{B}$ the set of dominant
integral weights for $\frak g_{B}$.
Denote by $\omega_{1}, \ldots , \omega_{4} \in P_{+}^{F}$ the
fundamental weights of $\frak g _{F}$, defined by
$\omega_{i}(\alpha_{j}^{\vee})= \delta _{ij}$ for all $i,j=1, \ldots ,4$,
and by $\bar{\omega}_{1}, \ldots , \bar{\omega}_{4} \in P_{+}^{B}$ the
fundamental weights of $\frak g _{B}$, defined by
$\bar{\omega}_{i}(\beta_{j}^{\vee})= \delta _{ij}$ for all $i,j=1, \ldots ,4$.

We shall also use the notation $V_{F}(\mu)$ (resp. $V_{B}(\mu)$)
for the irreducible highest weight ${\frak g}_{F}$-module
(resp. ${\frak g}_{B}$-module) with highest weight
$\mu \in {\frak h}^{*}$.

\section{Vertex operator algebras $L_{F}(n-\frac{7}{2},0)$ and
$L_{B}(n-\frac{7}{2},0)$, for $n \in \N$}
\label{sect.B-F.}

Let $\hat{\frak g}_{F}$, $\hat{\frak g}_{B}$, $\hat{\frak g}_{B}'$
and $\hat{\frak g}_{B}''$ be affine Lie algebras associated to
$\frak g _{F}$, $\frak g _{B}$, $\frak g _{B}'$ and $\frak g
_{B}''$, respectively. For $k \in \C$, denote by $N_{F}(k,0)$,
$N_{B}(k,0)$, $N_{B}'(k,0)$ and $N_{B}''(k,0)$ generalized Verma
modules associated to $\hat{\frak g}_{F}$, $\hat{\frak g}_{B}$,
$\hat{\frak g}_{B}'$ and $\hat{\frak g}_{B}''$ of level $k$, and by
$L_{F}(k,0)$, $L_{B}(k,0)$, $L_{B}'(k,0)$ and $L_{B}''(k,0)$
corresponding irreducible modules.

We shall also use the notation $M_{F}( \lambda)$ (resp. $M_{B}( \lambda)$)
for the Verma module for $\hat{\frak g}_{F}$
(resp. $\hat{\frak g}_{B}$) with highest
weight $\lambda\in \hat{{\frak h}}^{*}$, and $L_{F}( \lambda)$
(resp. $L_{B}( \lambda)$) for the irreducible
$\hat{\frak g}_{F}$-module (resp. $\hat{\frak g}_{B}$-module)
with highest weight $\lambda\in \hat{{\frak h}}^{*}$.

Since $\frak g _{B}$, $\frak g _{B}'$ and $\frak g _{B}''$ are Lie
subalgebras of $\frak g _{F}$, it follows that $\hat{\frak g}_{B}$,
$\hat{\frak g} _{B}'$ and $\hat{\frak g} _{B}''$ are Lie subalgebras
of $\hat{\frak g}_{F}$. Using the P-B-W theorem (which gives an
embedding of the universal enveloping algebra of a Lie subalgebra
into the universal enveloping algebra of a Lie algebra) we obtain
embeddings of generalized Verma modules $N_{B}(k,0)$, $N_{B}'(k,0)$
and $N_{B}''(k,0)$ into $N_{F}(k,0)$. Therefore, $N_{B}(k,0)$,
$N_{B}'(k,0)$ and $N_{B}''(k,0)$ are vertex subalgebras of
$N_{F}(k,0)$.

Furthermore, isomorphisms $\pi '$ and $\pi ''$ induce isomorphisms
of affine Lie algebras $\pi ' : \hat{\frak g} _{B} \to \hat{\frak g}
_{B}'$, $\pi '' : \hat{\frak g} _{B} \to \hat{\frak g} _{B}''$ and
isomorphisms of vertex operator algebras $\pi ' : N_{B}(k,0) \to
N_{B}'(k,0)$, $\pi '' : N_{B}(k,0) \to N_{B}''(k,0)$. Thus,
$N_{F}(k,0)$ contains three copies of $N_{B}(k,0)$ as vertex
subalgebras.

\subsection{Vertex operator algebra $L_{B}(n-\frac{7}{2},0)$, for $n \in \N$}

In this subsection we recall some facts about vertex operator
algebra $L_{B}(n-\frac{7}{2},0)$, for $n \in \N$. It was proved in
\cite[Theorem 11]{P} that the maximal $\hat{\frak g}_{B}$-submodule
of $N_{B}(n-\frac{7}{2},0)$ is generated by one singular vector:

\begin{prop} \label{p.3.1} The maximal $\hat{\frak g}_{B}$-submodule
of $N_{B}(n-\frac{7}{2},0)$ is $J_{B}(n-\frac{7}{2},0)=
U(\hat{\frak g}_{B})v_{B}$, where
\begin{eqnarray*}
&& v_{B}=\Big( -\frac{1}{4}e_{\epsilon_1}(-1)^2+ e_{\epsilon_1 - \epsilon_2}(-1)
e_{\epsilon_1 + \epsilon_2}(-1)+ e_{\epsilon_1 - \epsilon_3}(-1)
e_{\epsilon_1 + \epsilon_3}(-1) \\
&& \qquad \quad +e_{\epsilon_1 - \epsilon_4}(-1)
e_{\epsilon_1 + \epsilon_4}(-1) \Big)^n{\bf 1}
\end{eqnarray*}
is a singular vector in $N_{B}(n-\frac{7}{2},0)$.
\end{prop}

\begin{rem} The choice of root vectors for $\frak g _{B}$ in
(\ref{rel.kor.vekt.}) is slightly different from
the choice in \cite{P}, but formula
for the vector $v_{B}$ is the same in both bases.
\end{rem}
Thus
\[
L_{B}(n-\frac{7}{2},0) \cong \frac{N_{B}(n-\frac{7}{2},0)}
{U(\hat{\frak g}_{B})v_{B}}.
\]
\begin{coro} \label{c.3.0.1}
The associative algebra $A(L_{B}(n-\frac{7}{2},0))$ is
isomorphic to the algebra $U(\frak g_{B})/I_{B}$, where $I_{B}$
is the two-sided ideal of $U(\frak g_{B})$ generated by
\begin{eqnarray*}
u_{B}=\left(-\frac{1}{4}e_{\epsilon_1}^2+ e_{\epsilon_1 - \epsilon_2}
e_{\epsilon_1 + \epsilon_2}+ e_{\epsilon_1 - \epsilon_3}
e_{\epsilon_1 + \epsilon_3} +e_{\epsilon_1 - \epsilon_4}
e_{\epsilon_1 + \epsilon_4} \right)^n.
\end{eqnarray*}
\end{coro}

It follows from \cite[Theorem 24, Theorem 27]{P} that

\begin{prop} \label{p-klas-B}
\item[(1)] The set
\begin{eqnarray*}
\left\{ L_{B}(n-\frac{7}{2}, \mu) \ \vert \ \mu \in P_{+}^{B}, \ ( \mu , \epsilon_{1} ) \leq
n-\frac{1}{2} \right\}
\end{eqnarray*}
provides the complete list of irreducible
$L_{B}(n-\frac{7}{2},0)$-modules.
\item[(2)] Let $M$ be a $L_{B}(n-\frac{7}{2},0)$-module. Then
$M$ is completely reducible.
\end{prop}

In the special case $n=1$, the following was proved (\cite[Theorem
31, Theorem 33]{P}):

\begin{prop} \label{p-B-n=1}
\item[(1)] The set
\begin{eqnarray*}
&& \Big\{ L_{B}(-\frac{5}{2}, 0), L_{B}(-\frac{5}{2}, \bar{\omega} _{4}),
L_{B}(-\frac{5}{2}, -\frac{5}{2} \bar{\omega} _{1}),
L_{B}(-\frac{5}{2}, -\frac{7}{2} \bar{\omega} _{1}+ \bar{\omega} _{4}),
L_{B}(-\frac{5}{2}, -\frac{3}{2} \bar{\omega} _{2}), \\
&& L_{B}(-\frac{5}{2}, -\frac{5}{2} \bar{\omega} _{2}+ \bar{\omega} _{4}),
L_{B}(-\frac{5}{2}, -\frac{1}{2} \bar{\omega} _{3}),
L_{B}(-\frac{5}{2}, -\frac{3}{2} \bar{\omega} _{3}+ \bar{\omega} _{4}),
L_{B}(-\frac{5}{2}, \frac{1}{2} \bar{\omega} _{1}-
\frac{3}{2} \bar{\omega} _{2}), \\
&& L_{B}(-\frac{5}{2}, \frac{3}{2} \bar{\omega}_{1}- \frac{5}{2} \bar{\omega} _{2}
+ \bar{\omega} _{4}),
L_{B}(-\frac{5}{2}, -\frac{3}{2} \bar{\omega} _{1}- \frac{1}{2} \bar{\omega} _{3}),
L_{B}(-\frac{5}{2}, -\frac{1}{2} \bar{\omega} _{1}- \frac{3}{2} \bar{\omega} _{3}
+ \bar{\omega} _{4}), \\
&& L_{B}(-\frac{5}{2}, -\frac{1}{2} \bar{\omega} _{2}- \frac{1}{2} \bar{\omega} _{3}),
L_{B}(-\frac{5}{2}, \frac{1}{2} \bar{\omega} _{2}- \frac{3}{2} \bar{\omega} _{3}
+ \bar{\omega} _{4}),
L_{B}(-\frac{5}{2}, -\frac{1}{2} \bar{\omega} _{1}- \frac{1}{2} \bar{\omega} _{2}
- \frac{1}{2} \bar{\omega} _{3}), \\
&& L_{B}(-\frac{5}{2}, -\frac{3}{2} \bar{\omega} _{1}+ \frac{1}{2} \bar{\omega} _{2}
- \frac{3}{2} \bar{\omega} _{3}+ \bar{\omega} _{4}) \Big\}
\end{eqnarray*}
provides the complete list of irreducible weak
$L_{B}(-\frac{5}{2},0)$-modules from the category $\mathcal{O}$.
\item[(2)] Let $M$ be a weak $L_{B}(-\frac{5}{2},0)$-module from the category
$\mathcal{O}$. Then $M$ is completely reducible.
\end{prop}

\begin{rem} \label{r-B-serija} It follows from Propositions \ref{p-klas-B} and \ref{p-B-n=1}
that there are $16$ irreducible weak $L_{B}(-\frac{5}{2},0)$-modules from the category
$\mathcal{O}$, and that only irreducible $L_{B}(-\frac{5}{2},0)$-modules are
$L_{B}(-\frac{5}{2},0)$ and $L_{B}(-\frac{5}{2}, \bar{\omega} _{4})$.
\end{rem}

Using isomorphisms $\pi '$ and
$\pi ''$, and Proposition \ref{p.3.1} we obtain
\begin{prop} \label{prop.B.kopija1}
\item[(1)]
The maximal $\hat{\frak g}_{B}'$-submodule
of $N_{B}'(n-\frac{7}{2},0)$ is $J_{B}'(n-\frac{7}{2},0)=
U(\hat{\frak g}_{B}')v_{B}'$, where
\begin{eqnarray*}
&& v_{B}'=\Big( -\frac{1}{4}e_{\frac{1}{2}( \epsilon_{1} + \epsilon_{2}
+ \epsilon_{3} - \epsilon_{4})}(-1)^2+ e_{\epsilon_2 + \epsilon_3}(-1)
e_{\epsilon_1 - \epsilon_4}(-1)+ e_{\epsilon_2 - \epsilon_4}(-1)
e_{\epsilon_1 + \epsilon_3}(-1) \\
&& \qquad \quad +e_{\epsilon_3 - \epsilon_4}(-1)
e_{\epsilon_1 + \epsilon_2}(-1) \Big)^n{\bf 1}.
\end{eqnarray*}
\item[(2)]
The associative algebra $A(L_{B}'(n-\frac{7}{2},0))$ is
isomorphic to the algebra $U(\frak g_{B}')/I_{B}'$, where $I_{B}'$
is the two-sided ideal of $U(\frak g_{B}')$ generated by
\begin{eqnarray*}
u_{B}'=\left(-\frac{1}{4}e_{\frac{1}{2}( \epsilon_{1} + \epsilon_{2}
+ \epsilon_{3} - \epsilon_{4})}^2+ e_{\epsilon_2 + \epsilon_3}
e_{\epsilon_1 - \epsilon_4}+ e_{\epsilon_2 - \epsilon_4}
e_{\epsilon_1 + \epsilon_3} +e_{\epsilon_3 - \epsilon_4}
e_{\epsilon_1 + \epsilon_2} \right)^n.
\end{eqnarray*}
\end{prop}
{\bf Proof:} It can easily be checked that $\pi '(v_{B})=v_{B}'$,
which implies the claim of proposition. $\;\;\;\;\Box$
\begin{prop} \label{prop.B.kopija2}
\item[(1)]
The maximal $\hat{\frak g}_{B}''$-submodule
of $N_{B}''(n-\frac{7}{2},0)$ is $J_{B}''(n-\frac{7}{2},0)=
U(\hat{\frak g}_{B}'')v_{B}''$, where
\begin{eqnarray*}
&& v_{B}''=\Big( -\frac{1}{4}e_{\frac{1}{2}( \epsilon_{1} + \epsilon_{2}
+ \epsilon_{3} + \epsilon_{4})}(-1)^2+ e_{\epsilon_2 + \epsilon_3}(-1)
e_{\epsilon_1 + \epsilon_4}(-1)+ e_{\epsilon_2 + \epsilon_4}(-1)
e_{\epsilon_1 + \epsilon_3}(-1) \\
&& \qquad \quad +e_{\epsilon_3 + \epsilon_4}(-1)
e_{\epsilon_1 + \epsilon_2}(-1) \Big)^n{\bf 1}.
\end{eqnarray*}
\item[(2)]
The associative algebra $A(L_{B}''(n-\frac{7}{2},0))$ is
isomorphic to the algebra $U(\frak g_{B}'')/I_{B}''$, where $I_{B}''$
is the two-sided ideal of $U(\frak g_{B}'')$ generated by
\begin{eqnarray*}
u_{B}''=\left(-\frac{1}{4}e_{\frac{1}{2}( \epsilon_{1} + \epsilon_{2}
+ \epsilon_{3} + \epsilon_{4})}^2+ e_{\epsilon_2 + \epsilon_3}
e_{\epsilon_1 + \epsilon_4}+ e_{\epsilon_2 + \epsilon_4}
e_{\epsilon_1 + \epsilon_3} +e_{\epsilon_3 + \epsilon_4}
e_{\epsilon_1 + \epsilon_2} \right)^n.
\end{eqnarray*}
\end{prop}
{\bf Proof:} It can easily be checked that $\pi ''(v_{B})=v_{B}''$,
which implies the claim of proposition. $\;\;\;\;\Box$

We obtain
\begin{eqnarray*}
&& L_{B}'(n-\frac{7}{2},0) \cong \frac{N_{B}'(n-\frac{7}{2},0)}
{U(\hat{\frak g}_{B}')v_{B}'} \quad \mbox{and} \\
&& L_{B}''(n-\frac{7}{2},0) \cong \frac{N_{B}''(n-\frac{7}{2},0)}
{U(\hat{\frak g}_{B}'')v_{B}''}.
\end{eqnarray*}

\subsection{Vertex operator algebra $L_{F}(n-\frac{7}{2},0)$, for $n \in \N$}

In this subsection we show that the maximal $\hat{\frak g}_{F}$-submodule
of $N_{F}(n-\frac{7}{2},0)$, for $n \in \N$, is generated by one
singular vector. We need two lemmas to prove that.

Denote by $\lambda _{n}$ the weight $(n-\frac{7}{2})\Lambda_{0}$.
Then $N_{F}(n-\frac{7}{2},0)$ is a quotient of the Verma module
$M_{F}(\lambda _{n})$ and $L_{F}(n-\frac{7}{2},0) \cong L_{F}(\lambda _{n})$.

\begin{lem} \label{l.3.4}
The weight $\lambda _{n}=(n-\frac{7}{2})\Lambda_{0}$
is admissible for $\hat{\frak g}_{F}$ and
\[
\hat{\Pi}^{\vee}_{\lambda _{n}}= \{ (\delta - \epsilon_{1})^{\vee},
\alpha_{1}^{\vee},\alpha_{2}^{\vee}, \alpha_{3}^{\vee} ,
\alpha_{4}^{\vee} \}.
\]
\end{lem}
{\bf Proof:}
Clearly
\begin{eqnarray*}
&& \langle \lambda _{n} + \rho,\alpha _{i}^{\vee}\rangle = 1 \ \ \mbox{for } 1 \leq i \leq 4, \\
&& \langle \lambda _{n} + \rho,\alpha _{0}^{\vee} \rangle
=n-\frac{5}{2},
\end{eqnarray*}
which implies
\begin{eqnarray*}
\langle \lambda _{n} + \rho,(\delta - \epsilon_{1})^{\vee} \rangle =
\langle \lambda _{n} + \rho, 2 \alpha _{0}^{\vee}+2 \alpha _{1}^{\vee}+
2 \alpha _{2}^{\vee}+ \alpha _{3}^{\vee} \rangle = 2n.
\end{eqnarray*}
The claim of lemma now follows easily.
$\;\;\;\;\Box$

\begin{lem} \label{l.3.5}
Vector
\begin{eqnarray*}
&& v_{F}=\Big( -\frac{1}{4}e_{\epsilon_1}(-1)^2+ e_{\epsilon_1 - \epsilon_2}(-1)
e_{\epsilon_1 + \epsilon_2}(-1)+ e_{\epsilon_1 - \epsilon_3}(-1)
e_{\epsilon_1 + \epsilon_3}(-1) \\
&& \qquad \quad +e_{\epsilon_1 - \epsilon_4}(-1)
e_{\epsilon_1 + \epsilon_4}(-1) \Big)^n{\bf 1}.
\end{eqnarray*}
is a singular vector in $N_{F}(n-\frac{7}{2},0)$.
\end{lem}
{\bf Proof:} Clearly $v_{F} \in N_{B}(n-\frac{7}{2},0)$,
and Proposition \ref{p.3.1} implies that
\begin{eqnarray*}
& & e_{i}(0).v_{F}=0, \ 1 \leq i \leq 3 \\
& & f_{\theta}(1).v_{F}=0,
\end{eqnarray*}
since $e_{1}(0),e_{2}(0),e_{3}(0),f_{\theta}(1) \in
\hat{\frak g}_{B}$. It remains to show that $e_{4}(0).v_{F}=0$,
which follows immediately from the fact that $e_{4}(0)$ commutes
with all vectors $e_{\epsilon_1}(-1), e_{\epsilon_1 -
\epsilon_2}(-1), e_{\epsilon_1 + \epsilon_2}(-1), e_{\epsilon_1 -
\epsilon_3}(-1), e_{\epsilon_1 + \epsilon_3}(-1), e_{\epsilon_1 -
\epsilon_4}(-1), e_{\epsilon_1 + \epsilon_4}(-1)$. $\;\;\;\;\Box$

\begin{thm} \label{t.3.6} The maximal $\hat{\frak g}_{F}$-submodule
of $N_{F}(n-\frac{7}{2},0)$ is $J_{F}(n-\frac{7}{2},0)=
U(\hat{\frak g}_{F})v_{F}$, where
\begin{eqnarray*}
&& v_{F}=\Big( -\frac{1}{4}e_{\epsilon_1}(-1)^2+ e_{\epsilon_1 - \epsilon_2}(-1)
e_{\epsilon_1 + \epsilon_2}(-1)+ e_{\epsilon_1 - \epsilon_3}(-1)
e_{\epsilon_1 + \epsilon_3}(-1) \\
&& \qquad \quad +e_{\epsilon_1 - \epsilon_4}(-1)
e_{\epsilon_1 + \epsilon_4}(-1) \Big)^n{\bf 1}.
\end{eqnarray*}
\end{thm}
{\bf Proof:} It follows from Proposition \ref{t.KW1} and Lemma \ref{l.3.4}
that the maximal submodule of the Verma module $M_{F}(\lambda _{n})$ is generated by
five singular vectors with weights
\[
r_{\delta - \epsilon_1}.\lambda _{n}, \ r_{\alpha_{1}}.\lambda
_{n}, \ r_{\alpha_{2}}.\lambda _{n}, \ r_{\alpha_{3}}.\lambda _{n},
\ r_{\alpha_{4}}.\lambda _{n}.
\]
It follows from Lemma \ref{l.3.5} that $v_{F}$ is a singular vector
with weight
$\lambda _{n}-2n \delta +2n \epsilon_1=r_{\delta -
\epsilon_1}.\lambda _{n}$. Other singular vectors
have weights
\[
r_{\alpha_{i}}.\lambda _{n}=\lambda _{n}- \langle \lambda _{n}+
\rho, \alpha_{i}^{\vee} \rangle \alpha_{i}=\lambda _{n}
- \alpha_{i}, \ 1 \leq i \leq 4,
\]
so the images of these vectors under the projection of
$M_{F}(\lambda _{n})$ onto $N_{F}(n-\frac{7}{2},0)$ are 0.
Therefore, the maximal submodule of $N_{F}(n-\frac{7}{2},0)$ is
generated by the vector $v_{F}$, i.e. $J_{F}(n-\frac{7}{2},0)=
U(\hat{\frak g}_{F})v_{F}$. $\;\;\;\;\Box$

It follows that
\[
L_{F}(n-\frac{7}{2},0) \cong \frac{N_{F}(n-\frac{7}{2},0)}
{U(\hat{\frak g}_{F})v_{F}}.
\]

Using Theorem \ref{t.3.6} and Proposition \ref{p.1.5.5}
we can identify the associative algebra $A(L_{F}(n-\frac{7}{2},0))$:

\begin{prop} \label{p.3.7}
The associative algebra $A(L_{F}(n-\frac{7}{2},0))$ is
isomorphic to the algebra $U(\frak g_{F})/I_{F}$, where $I_{F}$
is the two-sided ideal of $U(\frak g_{F})$ generated by
\begin{eqnarray*}
u_{F}=\left(-\frac{1}{4}e_{\epsilon_1}^2+ e_{\epsilon_1 - \epsilon_2}
e_{\epsilon_1 + \epsilon_2}+ e_{\epsilon_1 - \epsilon_3}
e_{\epsilon_1 + \epsilon_3} +e_{\epsilon_1 - \epsilon_4}
e_{\epsilon_1 + \epsilon_4} \right)^n.
\end{eqnarray*}
\end{prop}

Since $v_{F}=v_{B}$, it follows from Theorem \ref{t.3.6} and Proposition \ref{p.3.1}
that the embedding of $N_{B}(n-\frac{7}{2},0)$ into $N_{F}(n-\frac{7}{2},0)$
induces the embedding of $L_{B}(n-\frac{7}{2},0)$ into $L_{F}(n-\frac{7}{2},0)$.
We get:

\begin{prop} \label{prop.ulaganje1}
$L_{B}(n-\frac{7}{2},0)$ is a vertex subalgebra of $L_{F}(n-\frac{7}{2},0)$.
\end{prop}
Denote by $R^{F}$ (resp. $R^{B}$) the $U(\frak g _{F})$-submodule
(resp. $U(\frak g _{B})$-submodule)
of $U(\frak g _{F})$ (resp. $U(\frak g _{B})$) generated by the
vector $u_{F}$ (resp. $u_{B}$) under the adjoint action. Let $R_{0}^{F}$
(resp. $R_{0}^{B}$) be the zero-weight
subspace of $R^{F}$ (resp. $R^{B}$). Denote by
${\mathcal P}_{0}^{F}$ and ${\mathcal P}_{0}^{B}$ the
corresponding sets of polynomials, defined in Subsection
\ref{subsec.3.3}. Since $u_{B}=u_{F}$, we have $R^{B} \subseteq R^{F}$
and
\begin{coro} \label{cor.ulaganje-pol1}
\[
{\mathcal P}_{0}^{B} \subseteq {\mathcal P}_{0}^{F}
\]
\end{coro}
Furthermore
\begin{prop} \label{prop.ulaganje2}
$L_{B}'(n-\frac{7}{2},0)$ and $L_{B}''(n-\frac{7}{2},0)$
are vertex subalgebras of $L_{F}(n-\frac{7}{2},0)$.
\end{prop}
{\bf Proof:} One can easily check that
\begin{eqnarray}
&& \frac{1}{(2n)!} f_{\frac{1}{2}( \epsilon_{1} - \epsilon_{2}
- \epsilon_{3} + \epsilon_{4})}(0)^{2n}. v_{F}= v_{B}' \quad
 \mbox{and} \nonumber  \\
&& \frac{1}{(2n)!} f_{\frac{1}{2}( \epsilon_{1} - \epsilon_{2}
- \epsilon_{3} - \epsilon_{4})}(0)^{2n}. v_{F}= v_{B}''.
\label{rel.ulaganje}
\end{eqnarray}
Theorem \ref{t.3.6} and Propositions \ref{prop.B.kopija1} and \ref{prop.B.kopija2}
now imply that the embedding of $N_{B}'(n-\frac{7}{2},0)$ into
$N_{F}(n-\frac{7}{2},0)$ induces the embedding of
$L_{B}'(n-\frac{7}{2},0)$ into $L_{F}(n-\frac{7}{2},0)$,
and that the embedding of $N_{B}''(n-\frac{7}{2},0)$ into
$N_{F}(n-\frac{7}{2},0)$ induces the embedding of
$L_{B}''(n-\frac{7}{2},0)$ into $L_{F}(n-\frac{7}{2},0)$.
$\;\;\;\;\Box$

We conclude that $L_{F}(n-\frac{7}{2},0)$ contains three copies of
$L_{B}(n-\frac{7}{2},0)$ as vertex subalgebras.

Denote by $R^{B'}$ (resp. $R^{B''}$) the $U(\frak g _{B}')$-submodule
(resp. $U(\frak g _{B}'')$-submodule)
of $U(\frak g _{B}')$ (resp. $U(\frak g _{B}'')$) generated by the
vector $u_{B}'$ (resp. $u_{B}''$) under the adjoint action. Let $R_{0}^{B'}$
(resp. $R_{0}^{B''}$) be the zero-weight
subspace of $R^{B'}$ (resp. $R^{B''}$). Denote by
${\mathcal P}_{0}^{B'}$ and ${\mathcal P}_{0}^{B''}$ the
corresponding sets of polynomials, defined in Subsection
\ref{subsec.3.3}.
It follows from relations (\ref{rel.ulaganje}) that
\begin{eqnarray*}
&& \frac{1}{(2n)!} (f_{\frac{1}{2}( \epsilon_{1} - \epsilon_{2}
- \epsilon_{3} + \epsilon_{4})}^{2n}) _{L} u_{F}= u_{B}' \quad
 \mbox{and}  \\
&& \frac{1}{(2n)!} (f_{\frac{1}{2}( \epsilon_{1} - \epsilon_{2}
- \epsilon_{3} - \epsilon_{4})}^{2n}) _{L} u_{F}= u_{B}'',
\end{eqnarray*}
which implies that $R^{B'} \subseteq R^{F}$,
$R^{B''} \subseteq R^{F}$ and

\begin{coro} \label{cor.ulaganje-pol2}
\[
{\mathcal P}_{0}^{B'} \subseteq {\mathcal P}_{0}^{F} \quad \mbox{and}
\quad {\mathcal P}_{0}^{B''} \subseteq {\mathcal P}_{0}^{F}
\]
\end{coro}

\section{Vertex operator algebras $L_{B}(-\frac{5}{2},0)$ and
$L_{F}(-\frac{5}{2},0)$ }
\label{sect.n=1}

In this section we study the case $n=1$, vertex operator
algebras $L_{B}(-\frac{5}{2},0)$ and $L_{F}(-\frac{5}{2},0)$.
We show that $L_{B}(-\frac{5}{2},0)$ is a vertex operator subalgebra
of $L_{F}(-\frac{5}{2},0)$, i.e. that these
vertex operator algebras have the same conformal vector.
We use the following lemma:
\begin{lem} Relation
\begin{eqnarray} \label{rel.dokaz.Vir}
&& 7 \sum_{\alpha \in \Delta _{B}^{+} \atop (\alpha,\alpha)=1 }( e_{\alpha}(-1)f_{\alpha}(-1)
{\bf 1}+f_{\alpha}(-1)e_{\alpha}(-1) {\bf 1}) \nonumber \\
&& =4\sum_{\alpha \in \Delta _{B}^{+} \atop (\alpha,\alpha)=2 }( e_{\alpha}(-1)
f_{\alpha}(-1){\bf 1} +f_{\alpha}(-1)e_{\alpha}(-1) {\bf 1})
+ \sum_{\alpha \in \Delta _{B}^{+} \atop (\alpha,\alpha)=1 }
h_{\alpha}(-1)^{2}{\bf 1} \qquad
\end{eqnarray}
holds in $L_{B}(-\frac{5}{2},0)$.
\end{lem}
{\bf Proof:} In the case $n=1$, Proposition \ref{p.3.1} implies that
relation
\begin{eqnarray*}
&& \Big( 2f_{\epsilon_1}(0)^2 + f_{\epsilon_2}(0)^2
f_{\epsilon_1 - \epsilon_2}(0)^2 + f_{\epsilon_3}(0)^2
f_{\epsilon_1 - \epsilon_3}(0)^2 + f_{\epsilon_4}(0)^2
f_{\epsilon_1 - \epsilon_4}(0)^2 \Big). \\
&& \Big( -\frac{1}{4}e_{\epsilon_1}(-1)^2 {\bf 1}+ e_{\epsilon_1 - \epsilon_2}(-1)
e_{\epsilon_1 + \epsilon_2}(-1) {\bf 1}+ e_{\epsilon_1 - \epsilon_3}(-1)
e_{\epsilon_1 + \epsilon_3}(-1) {\bf 1} \\
&& \quad +e_{\epsilon_1 - \epsilon_4}(-1)
e_{\epsilon_1 + \epsilon_4}(-1){\bf 1} \Big) =0 \nonumber
\end{eqnarray*}
holds in $L_{B}(-\frac{5}{2},0)$. From this relation we get
the claim of lemma. $\;\;\;\;\Box$

\begin{thm} \label{prop.jedn.Vir}
Denote by $\omega _{F}$ the conformal vector for vertex operator
algebra $L_{F}(-\frac{5}{2},0)$, and by $\omega _{B}$ the
conformal vector for $L_{B}(-\frac{5}{2},0)$. Then
\begin{eqnarray*}
\omega _{F}=\omega _{B}.
\end{eqnarray*}
\end{thm}
{\bf Proof:} Sets $\{ \frac{1}{2} h_{\alpha} \ \vert \ \alpha \in \Delta
_{B}^{+}, (\alpha,\alpha)=1 \}$, $\{ \frac{1}{2} h_{\alpha} \ \vert \ \alpha \in \Delta
_{B}'^{+}, (\alpha,\alpha)=1 \}$ and $\{ \frac{1}{2} h_{\alpha} \ \vert \ \alpha \in \Delta
_{B}''^{+}, (\alpha,\alpha)=1 \}$ are three
orthonormal bases of ${\frak h}$ with respect to the form $(\cdot, \cdot)$.
Clearly,
\begin{eqnarray} \label{rel.5.0.1}
\sum_{\alpha \in \Delta _{B}^{+} \atop (\alpha,\alpha)=1 }
h_{\alpha}(-1)^{2}{\bf 1}=\sum_{\alpha \in \Delta _{B}'^{+} \atop (\alpha,\alpha)=1 }
h_{\alpha}(-1)^{2}{\bf 1}=\sum_{\alpha \in \Delta _{B}''^{+} \atop (\alpha,\alpha)=1 }
h_{\alpha}(-1)^{2}{\bf 1}.
\end{eqnarray}
Applying isomorphisms $\pi '$ and $\pi ''$ to relation
(\ref{rel.dokaz.Vir}) we obtain that relation
\begin{eqnarray} \label{rel.dokaz.Vir2}
&& 7 \sum_{\alpha \in \Delta _{B}'^{+} \atop (\alpha,\alpha)=1 }( e_{\alpha}(-1)f_{\alpha}(-1)
{\bf 1}+f_{\alpha}(-1)e_{\alpha}(-1) {\bf 1}) \nonumber \\
&& =4\sum_{\alpha \in \Delta _{B}^{+} \atop (\alpha,\alpha)=2 }( e_{\alpha}(-1)
f_{\alpha}(-1){\bf 1} +f_{\alpha}(-1)e_{\alpha}(-1) {\bf 1})
+ \sum_{\alpha \in \Delta _{B}'^{+} \atop (\alpha,\alpha)=1 }
h_{\alpha}(-1)^{2}{\bf 1} \qquad
\end{eqnarray}
holds in the vertex subalgebra $L_{B}'(-\frac{5}{2},0)$ of
$L_{F}(-\frac{5}{2},0)$, and that relation
\begin{eqnarray} \label{rel.dokaz.Vir3}
&& 7 \sum_{\alpha \in \Delta _{B}''^{+} \atop (\alpha,\alpha)=1 }( e_{\alpha}(-1)f_{\alpha}(-1)
{\bf 1}+f_{\alpha}(-1)e_{\alpha}(-1) {\bf 1}) \nonumber \\
&& =4\sum_{\alpha \in \Delta _{B}^{+} \atop (\alpha,\alpha)=2 }( e_{\alpha}(-1)
f_{\alpha}(-1){\bf 1} +f_{\alpha}(-1)e_{\alpha}(-1) {\bf 1})
+ \sum_{\alpha \in \Delta _{B}''^{+} \atop (\alpha,\alpha)=1 }
h_{\alpha}(-1)^{2}{\bf 1} \qquad
\end{eqnarray}
holds in the vertex subalgebra $L_{B}''(-\frac{5}{2},0)$ of
$L_{F}(-\frac{5}{2},0)$. Using (\ref{rel.5.0.1}) we obtain that
\begin{eqnarray} \label{rel.dokaz.Vir4}
&& \sum_{\alpha \in \Delta _{B}^{+} \atop (\alpha,\alpha)=1 }( e_{\alpha}(-1)f_{\alpha}(-1)
{\bf 1}+f_{\alpha}(-1)e_{\alpha}(-1) {\bf 1}) \nonumber \\
&& =\sum_{\alpha \in \Delta _{B}'^{+} \atop (\alpha,\alpha)=1 }( e_{\alpha}(-1)f_{\alpha}(-1)
{\bf 1}+f_{\alpha}(-1)e_{\alpha}(-1) {\bf 1}) \nonumber \\
&&=\sum_{\alpha \in \Delta _{B}''^{+} \atop (\alpha,\alpha)=1 }( e_{\alpha}(-1)f_{\alpha}(-1)
{\bf 1}+f_{\alpha}(-1)e_{\alpha}(-1) {\bf 1}) \nonumber \\
&& = \frac{4}{7} \sum_{\alpha \in \Delta _{B}^{+} \atop (\alpha,\alpha)=2 }( e_{\alpha}(-1)
f_{\alpha}(-1){\bf 1} +f_{\alpha}(-1)e_{\alpha}(-1) {\bf 1})
+ \frac{1}{7} \sum_{\alpha \in \Delta _{B}^{+} \atop (\alpha,\alpha)=1 }
h_{\alpha}(-1)^{2}{\bf 1} \qquad
\end{eqnarray}
holds in $L_{F}(-\frac{5}{2},0)$.
It follows from
relation (\ref{rel.Virasoro}) that
\begin{eqnarray*}
&& \omega _{F}=\frac{1}{13} \Bigg(\frac{1}{4} \sum_{\alpha \in \Delta _{B}^{+} \atop (\alpha,\alpha)=1 }
h_{\alpha}(-1)^{2}{\bf 1} + \sum_{\alpha \in \Delta _{F}^{+}}
\frac{(\alpha, \alpha)}{2}(e_{\alpha}(-1)f_{\alpha}(-1) {\bf 1}+
f_{\alpha}(-1)e_{\alpha}(-1) {\bf 1}) \Bigg) \\
&& = \frac{1}{13} \Bigg( \frac{1}{4} \sum_{\alpha \in \Delta _{B}^{+} \atop (\alpha,\alpha)=1 }
h_{\alpha}(-1)^{2}{\bf 1}+
\sum_{\alpha \in \Delta _{B}^{+} \atop (\alpha,\alpha)=2 }
(e_{\alpha}(-1)f_{\alpha}(-1) {\bf 1}+ f_{\alpha}(-1)e_{\alpha}(-1) {\bf 1}) \\
&&  \qquad + \frac{1}{2} \sum_{\alpha \in \Delta _{B}^{+} \atop (\alpha,\alpha)=1 }( e_{\alpha}(-1)f_{\alpha}(-1)
{\bf 1}+f_{\alpha}(-1)e_{\alpha}(-1) {\bf 1}) \\
&& \qquad + \frac{1}{2} \sum_{\alpha \in \Delta _{B}'^{+} \atop (\alpha,\alpha)=1 }( e_{\alpha}(-1)f_{\alpha}(-1)
{\bf 1}+f_{\alpha}(-1)e_{\alpha}(-1) {\bf 1}) \\
&& \qquad + \frac{1}{2} \sum_{\alpha \in \Delta _{B}''^{+} \atop (\alpha,\alpha)=1 }( e_{\alpha}(-1)f_{\alpha}(-1)
{\bf 1}+f_{\alpha}(-1)e_{\alpha}(-1) {\bf 1}) \Bigg).
\end{eqnarray*}
Using relation (\ref{rel.dokaz.Vir4}), we obtain
\begin{eqnarray*}
&& \omega _{F}=
\frac{1}{28} \sum_{\alpha \in \Delta _{B}^{+} \atop (\alpha,\alpha)=1 }
h_{\alpha}(-1)^{2}{\bf 1}+
\frac{1}{7} \sum_{\alpha \in \Delta _{B}^{+} \atop (\alpha,\alpha)=2 }
(e_{\alpha}(-1)f_{\alpha}(-1) {\bf 1}+ f_{\alpha}(-1)e_{\alpha}(-1) {\bf 1}) \\
&& = \frac{1}{9} \Bigg( \frac{1}{4} \sum_{\alpha \in \Delta _{B}^{+} \atop (\alpha,\alpha)=1 }
h_{\alpha}(-1)^{2}{\bf 1}+
\sum_{\alpha \in \Delta _{B}^{+} \atop (\alpha,\alpha)=2 }
(e_{\alpha}(-1)f_{\alpha}(-1) {\bf 1}+ f_{\alpha}(-1)e_{\alpha}(-1) {\bf 1}) \\
&&  \qquad + \frac{1}{2} \sum_{\alpha \in \Delta _{B}^{+} \atop (\alpha,\alpha)=1 }( e_{\alpha}(-1)f_{\alpha}(-1)
{\bf 1}+f_{\alpha}(-1)e_{\alpha}(-1) {\bf 1}) \Bigg) =
\omega _{B}. \;\;\;\;\Box
\end{eqnarray*}

Thus, $L_{B}(-\frac{5}{2},0)$ is a vertex operator subalgebra
of $L_{F}(-\frac{5}{2},0)$.

\begin{rem} Denote by $\omega _{B}'$ the conformal vector for vertex operator
algebra $L_{B}'(-\frac{5}{2},0)$, and by $\omega _{B}''$ the
conformal vector for $L_{B}''(-\frac{5}{2},0)$. Relation
(\ref{rel.dokaz.Vir4}) implies that
\[
\omega _{B}'=\omega _{B}''=\omega _{B}=\omega _{F}
\]
in $L_{F}(-\frac{5}{2},0)$.
Thus, $L_{B}'(-\frac{5}{2},0)$ and
$L_{B}''(-\frac{5}{2},0)$ are vertex operator subalgebras of
$L_{F}(-\frac{5}{2},0)$.
\end{rem}

In the next theorem we determine the decomposition of
$L_{F}(-\frac{5}{2},0)$ into a direct sum of irreducible
$L_{B}(-\frac{5}{2},0)$-modules.

\begin{thm} \label{thm.extension}
\begin{eqnarray*}
L_{F}(-\frac{5}{2},0) \cong L_{B}(-\frac{5}{2},0) \oplus
L_{B}(-\frac{5}{2}, \bar{\omega} _{4})
\end{eqnarray*}
\end{thm}
{\bf Proof:} It follows from Theorem \ref{prop.jedn.Vir} that $L_{F}(-\frac{5}{2},0)$ is a
$L_{B}(-\frac{5}{2},0)$-module, and Proposition \ref{p-klas-B}
implies that it is a direct sum of copies
of irreducible $L_{B}(-\frac{5}{2},0)$-modules $L_{B}(-\frac{5}{2},0)$
and $L_{B}(-\frac{5}{2}, \bar{\omega} _{4})$.
Clearly, ${\bf 1}$ and $e_{\frac{1}{2}( \epsilon_{1} + \epsilon_{2}
+ \epsilon_{3} + \epsilon_{4})}(-1){\bf 1}$ are singular vectors
for $\hat{\frak g}_{B}$ in $L_{F}(-\frac{5}{2},0)$ which generate
the following $L_{B}(-\frac{5}{2},0)$-modules:
\begin{eqnarray*}
&& U(\hat{\frak g}_{B}). {\bf 1} \cong L_{B}(-\frac{5}{2},0) \quad \mbox{and} \\
&& U(\hat{\frak g}_{B}). e_{\frac{1}{2}( \epsilon_{1} + \epsilon_{2}
+ \epsilon_{3} + \epsilon_{4})}(-1){\bf 1} \cong L_{B}(-\frac{5}{2},
\bar{\omega} _{4}).
\end{eqnarray*}
It follows from relation (\ref{rel.lowest.conf.w})
that the lowest conformal weight of $L_{B}(-\frac{5}{2},0)$-module
$L_{B}(-\frac{5}{2},0)$ is $0$ and of $L_{B}(-\frac{5}{2},
\bar{\omega} _{4})$ is $1$. Theorem \ref{prop.jedn.Vir}
now implies that ${\bf 1}$ and $e_{\frac{1}{2}( \epsilon_{1} + \epsilon_{2}
+ \epsilon_{3} + \epsilon_{4})}(-1){\bf 1}$ are only singular vectors
for $\hat{\frak g}_{B}$ in $L_{F}(-\frac{5}{2},0)$, which implies
the claim of theorem. $\;\;\;\;\Box$

\begin{rem}
\item[(a)] It follows from Theorem \ref{thm.extension}
that the extension of vertex operator algebra
$L_{B}(-\frac{5}{2},0)$ by $L_{B}(-\frac{5}{2}, \bar{\omega} _{4})$,
its only irreducible module other then itself,
is a vertex operator algebra,
which is isomorphic to $L_{F}(-\frac{5}{2},0)$.

\item[(b)] Theorem \ref{thm.extension} also implies that
$\hat{\frak g}_{F}$-module $L_{F}(-\frac{5}{2},0)$, considered
as a module for Lie subalgebra $\hat{\frak g}_{B}$ of
$\hat{\frak g}_{F}$, decomposes into the finite direct sum
of $\hat{\frak g}_{B}$-modules. In the case of positive integer
levels, such cases are in physics called conformal embeddings of affine Lie
algebras and were studied in \cite{AGO}, \cite{BB}, \cite{SW}.
\end{rem}

\section{Weak $L_{F}(-\frac{5}{2},0)$-modules
from category $\mathcal{O}$} \label{sect.kateg.O}

In this section we study the category of weak
$L_{F}(-\frac{5}{2},0)$-modules that are in category $\mathcal{O}$
as $\hat{\frak g}_{F}$-modules. To obtain the classification of
irreducible objects in that category, we use methods from \cite{MP},
\cite{A2}, \cite{AM} (presented in Subsection \ref{subsec.3.3}) and
the fact that $L_{F}(-\frac{5}{2},0)$ contains three copies of
$L_{B}(-\frac{5}{2},0)$ as vertex operator subalgebras. It is proved
in \cite[Lemma 18]{P} that

\begin{lem} \label{l.polinomiB}
Let
\begin{eqnarray*}
&& p_{1}(h)=h_{\epsilon_1 - \epsilon_{2}}(h_{\epsilon_1 +
\epsilon_{2}}+\frac{5}{2}), \\
&& p_{2}(h)=h_{\epsilon_2 - \epsilon_{3}}(h_{\epsilon_2 +
\epsilon_{3}}+\frac{3}{2}), \\
&& p_{3}(h)=h_{\epsilon_3 - \epsilon_{4}}(h_{\epsilon_3 +
\epsilon_{4}}+\frac{1}{2}), \\
&& p_{4}(h)=h_{\epsilon _4}(h_{\epsilon _4}-1).
\end{eqnarray*}
Then $p_{1}, p_{2}, p_{3}, p_{4} \in {\mathcal P}_{0}^{B}$.
\end{lem}

Corollary \ref{cor.ulaganje-pol1} now implies that
$p_{1}, p_{2}, p_{3}, p_{4} \in {\mathcal P}_{0}^{F}$.

\begin{lem} \label{l.polinomiF}
Let
\begin{eqnarray*}
&& p_{5}(h)=h_{\epsilon_3 - \epsilon_{4}}(h_{\epsilon_1 +
\epsilon_{2}}+\frac{5}{2}), \\
&& p_{6}(h)=h_{\epsilon_2 - \epsilon_{3}}(h_{\epsilon_1 +
\epsilon_{4}}+\frac{3}{2}), \\
&& p_{7}(h)=h_{\epsilon_3 + \epsilon_{4}}(h_{\epsilon_1 -
\epsilon_{2}}+\frac{1}{2}), \\
&& p_{8}(h)=h_{\frac{1}{2}(\epsilon_1 - \epsilon_{2}-
\epsilon_{3} - \epsilon_{4})}(h_{\frac{1}{2}(\epsilon_1 - \epsilon_{2}-
\epsilon_{3} - \epsilon_{4})}-1), \\
&& p_{9}(h)=h_{\epsilon_3 + \epsilon_{4}}(h_{\epsilon_1 +
\epsilon_{2}}+\frac{5}{2}), \\
&& p_{10}(h)=h_{\epsilon_2 - \epsilon_{3}}(h_{\epsilon_1 -
\epsilon_{4}}+\frac{3}{2}), \\
&& p_{11}(h)=h_{\epsilon_3 - \epsilon_{4}}(h_{\epsilon_1 -
\epsilon_{2}}+\frac{1}{2}), \\
&& p_{12}(h)=h_{\frac{1}{2}(\epsilon_1 - \epsilon_{2}-
\epsilon_{3} + \epsilon_{4})}(h_{\frac{1}{2}(\epsilon_1 - \epsilon_{2}-
\epsilon_{3} + \epsilon_{4})}-1).
\end{eqnarray*}
Then $p_{5}, \ldots , p_{12} \in {\mathcal P}_{0}^{F}$.
\end{lem}
{\bf Proof:} Using isomorphisms $\pi '$ and
$\pi ''$, and Lemma \ref{l.polinomiB} we obtain
$p_{5}, \ldots , p_{8} \in {\mathcal P}_{0}^{B'}$
and $p_{9}, \ldots , p_{12} \in {\mathcal P}_{0}^{B''}$.
Corollary \ref{cor.ulaganje-pol2} now implies that $p_{5}, \ldots , p_{12} \in {\mathcal P}_{0}^{F}$.
$\;\;\;\;\Box$

\begin{prop}  The set
\begin{eqnarray*}
\left\{ V_{F}(0), V_{F}(-\frac{3}{2} \omega _{1}), V_{F}(-\frac{1}{2} \omega _{1}
-\frac{1}{2} \omega _{2}), V_{F}(-\frac{3}{2} \omega _{2}
+ \omega _{3})  \right\}
\end{eqnarray*}
provides the complete list of irreducible $A(L_{F}(-\frac{5}{2},0))$-modules
from the category $\mathcal{O}$.
\end{prop}
{\bf Proof:} Corollary \ref{c.1.7.2} implies that highest weights
$\mu \in {\frak h}^{*}$ of irreducible $A(L_{F}(-\frac{5}{2},0))$-modules
from the category $\mathcal{O}$ are in one-to-one correspondence
with the solutions of the polynomial equations $p(\mu)=0$
for all $p \in {\mathcal P}_{0}^{F}$. Since $R^{F}$ is the
$U(\frak g _{F})$-submodule of $U(\frak g _{F})$  generated by the
vector $u_{F}$ under the adjoint action, it is isomorphic to
$V_{F}(2 \epsilon _{1})$, i.e. $V_{F}(2 \omega _{4})$.
One can easily check that $\dim R_{0}^{F}=12$, which
implies that $\dim {\mathcal P}_{0}^{F}=12$. Since
the polynomials $p_{1}, \ldots , p_{12}$ from
Lemmas \ref{l.polinomiB} and \ref{l.polinomiF} are
linearly independent elements in ${\mathcal P}_{0}^{F}$, they form a basis for
${\mathcal P}_{0}^{F}$. The claim of proposition now easily
follows by solving the system $p_{i}(\mu)=0$, for $i=1, \ldots ,12$.
$\;\;\;\;\Box$

It follows from Zhu's theory that:

\begin{thm} \label{c.2.6.3}
The set
\begin{eqnarray*}
\left\{ L_{F}(-\frac{5}{2}, 0), L_{F}(-\frac{5}{2}, -\frac{3}{2} \omega _{1}),
L_{F}(-\frac{5}{2}, -\frac{1}{2} \omega _{1}
-\frac{1}{2} \omega _{2}), L_{F}(-\frac{5}{2}, -\frac{3}{2} \omega _{2}
+ \omega _{3})  \right\}
\end{eqnarray*}
provides the complete list of irreducible weak $L_{F}(-\frac{5}{2},0)$-modules
from the category $\mathcal{O}$.
\end{thm}

Denote by $\lambda ^{1}= -\frac{5}{2} \Lambda _{0}$,
$\lambda ^{2}= -\frac{5}{2} \Lambda _{0}-\frac{3}{2} \omega _{1}$,
$\lambda ^{3}= -\frac{5}{2} \Lambda _{0}-\frac{1}{2} \omega _{1}
-\frac{1}{2} \omega _{2}$ and $\lambda ^{4}= -\frac{5}{2} \Lambda _{0}
-\frac{3}{2} \omega _{2} + \omega _{3}$ the highest weights of
irreducible weak $L_{F}(-\frac{5}{2},0)$-modules
from the category $\mathcal{O}$. The following lemma is
crucial for proving complete reducibility of weak
$L_{F}(-\frac{5}{2},0)$-modules from the category
$\mathcal{O}$.

\begin{prop} \label{prop.dopustivost}
Weights $\lambda ^{i}$, $i=1,2,3,4$ are admissible for $\hat{\frak
g}_{F}$.
\end{prop}
{\bf Proof:} It is already proved in Lemma \ref{l.3.4} that the
weight $\lambda ^{1}= \lambda _{1}=-\frac{5}{2} \Lambda _{0}$ is admissible.
Similarly, one can easily check that weights
$\lambda ^{2}$, $\lambda ^{3}$ and $\lambda ^{4}$ are admissible and that
\begin{eqnarray*}
&& \hat{\Pi}^{\vee}_{\lambda ^{2}}= \{ (\delta - \epsilon_{1}-\epsilon_{3})^{\vee},
\alpha_{2}^{\vee},\alpha_{3}^{\vee}, \alpha_{4}^{\vee} ,
\epsilon_{2}^{\vee} \}, \\
&& \hat{\Pi}^{\vee}_{\lambda ^{3}}= \{ \alpha_{0}^{\vee},
(\epsilon_{2}-\epsilon_{4})^{\vee}, \alpha_{3}^{\vee},\alpha_{4}^{\vee},
\epsilon_{3}^{\vee} \}, \\
&& \hat{\Pi}^{\vee}_{\lambda ^{4}}= \{ \alpha_{0}^{\vee}, \alpha_{1}^{\vee},
\epsilon_{3}^{\vee}, \alpha_{4}^{\vee}, \alpha_{3}^{\vee} \}.
\;\;\;\;\Box
\end{eqnarray*}
We obtain:
\begin{thm}
Let $M$ be a weak $L_{F}(-\frac{5}{2},0)$-module from the category
$\mathcal{O}$. Then $M$ is completely reducible.
\end{thm}
{\bf Proof:} Let $L_{F}(\lambda)$ be some irreducible subquotient
of $M$. Then $L_{F}(\lambda)$ is an irreducible
weak $L_{F}(-\frac{5}{2},0)$-module
and Theorem \ref{c.2.6.3} implies that
$\lambda=\lambda ^{i}$, for some $i \in \{1,2,3,4\}$.
It follows from Proposition \ref{prop.dopustivost} that $\lambda$ is admissible.
Proposition \ref{t.KW2} now implies that $M$ is completely reducible.
$\;\;\;\;\Box$

\section{Decompositions of irreducible weak $L_{F}(-\frac{5}{2},0)$-modules
from category $\mathcal{O}$ into direct sums of
irreducible weak $L_{B}(-\frac{5}{2},0)$-modules} \label{sect.decomp}

In this section we determine decompositions of
irreducible weak $L_{F}(-\frac{5}{2},0)$-modules
from category $\mathcal{O}$,
classified in Theorem \ref{c.2.6.3},
into direct sums of irreducible
weak $L_{B}(-\frac{5}{2},0)$-modules, and show that these
direct sums are finite. We use the fact that these
vertex operator algebras have the same conformal vector.

Using relation (\ref{rel.lowest.conf.w}), we can determine
the lowest conformal weights of
irreducible weak $L_{B}(-\frac{5}{2},0)$-modules from
category $\mathcal{O}$, listed in Proposition \ref{p-B-n=1}:

\begin{lem} \label{lem.lowest.conf.w}
\item[(1)] The lowest conformal weight of weak
$L_{B}(-\frac{5}{2},0)$-modules $L_{B}(-\frac{5}{2}, -\frac{5}{2} \bar{\omega}
_{1})$, $L_{B}(-\frac{5}{2}, -\frac{3}{2} \bar{\omega} _{3}+ \bar{\omega}
_{4})$, $L_{B}(-\frac{5}{2}, \frac{1}{2} \bar{\omega} _{1}-
\frac{3}{2} \bar{\omega} _{2})$, $L_{B}(-\frac{5}{2}, -\frac{1}{2} \bar{\omega} _{2}-
\frac{1}{2} \bar{\omega} _{3})$ is $-\frac{5}{4}$.
\item[(2)] The lowest conformal weight of
$L_{B}(-\frac{5}{2}, -\frac{7}{2} \bar{\omega} _{1}+ \bar{\omega}
_{4})$, $L_{B}(-\frac{5}{2}, -\frac{1}{2} \bar{\omega} _{3})$,
$L_{B}(-\frac{5}{2}, \frac{3}{2} \bar{\omega}_{1}- \frac{5}{2} \bar{\omega} _{2}
+ \bar{\omega} _{4})$ and $L_{B}(-\frac{5}{2}, \frac{1}{2}
\bar{\omega} _{2}- \frac{3}{2} \bar{\omega} _{3}
+ \bar{\omega} _{4})$ is $-\frac{3}{4}$.
\item[(3)] The lowest conformal weight of
$L_{B}(-\frac{5}{2}, -\frac{3}{2} \bar{\omega}
_{2})$, $L_{B}(-\frac{5}{2}, -\frac{5}{2} \bar{\omega} _{2}+ \bar{\omega}
_{4})$, $L_{B}(-\frac{5}{2}, -\frac{1}{2} \bar{\omega} _{1}
- \frac{1}{2} \bar{\omega} _{2} - \frac{1}{2} \bar{\omega} _{3})$,
$L_{B}(-\frac{5}{2}, -\frac{3}{2} \bar{\omega} _{1}+ \frac{1}{2} \bar{\omega} _{2}
- \frac{3}{2} \bar{\omega} _{3}+ \bar{\omega} _{4})$, $L_{B}(-\frac{5}{2}, -\frac{3}{2} \bar{\omega} _{1}
- \frac{1}{2} \bar{\omega} _{3})$ and
$L_{B}(-\frac{5}{2}, -\frac{1}{2} \bar{\omega} _{1}-
\frac{3}{2} \bar{\omega} _{3} + \bar{\omega} _{4})$ is
$-\frac{3}{2}$.
\item[(4)] The lowest conformal weight of $L_{B}(-\frac{5}{2},0)$ is
$0$ and of $L_{B}(-\frac{5}{2}, \bar{\omega} _{4})$ is $1$.
\end{lem}

Denote by $v^{(2)}$ the highest weight vector of $L_{F}(-\frac{5}{2},
-\frac{3}{2} \omega _{1})$, by $v^{(3)}$ the highest weight
vector of $L_{F}(-\frac{5}{2}, -\frac{1}{2} \omega _{1}
-\frac{1}{2} \omega _{2})$ and by $v^{(4)}$ the highest weight
vector of $L_{F}(-\frac{5}{2}, -\frac{3}{2} \omega _{2}
+ \omega _{3})$. Then

\begin{lem} \label{lem.7.2}
\item[(1)] Relation
\[
f_{\frac{1}{2}( \epsilon_{1} - \epsilon_{2}
+ \epsilon_{3} - \epsilon_{4})}(0)v^{(2)}=0
\]
holds in $L_{F}(-\frac{5}{2}, -\frac{3}{2} \omega _{1})$.
\item[(2)] Relation
\[
f_{\frac{1}{2}( \epsilon_{1} - \epsilon_{2}
- \epsilon_{3} + \epsilon_{4})}(0)v^{(3)}=0
\]
holds in $L_{F}(-\frac{5}{2}, -\frac{1}{2} \omega _{1}
-\frac{1}{2} \omega _{2})$.
\item[(3)] Relation
\[
f_{\frac{1}{2}( \epsilon_{1} - \epsilon_{2}
- \epsilon_{3} - \epsilon_{4})}(0)v^{(4)}=0
\]
holds in $L_{F}(-\frac{5}{2}, -\frac{3}{2} \omega _{2}
+ \omega _{3})$.
\end{lem}
{\bf Proof:} We shall prove relation (1). Relations (2) and (3) can
be proved similarly. It is proved in Proposition \ref{prop.dopustivost}
that the weight $\lambda ^{2}= -\frac{5}{2} \Lambda _{0}-\frac{3}{2} \omega _{1}$
is admissible and that
\[
\hat{\Pi}^{\vee}_{\lambda ^{2}}= \{ (\delta - \epsilon_{1}-\epsilon_{3})^{\vee},
\alpha_{2}^{\vee},\alpha_{3}^{\vee}, \alpha_{4}^{\vee} ,
\epsilon_{2}^{\vee} \}.
\]
Proposition \ref{t.KW1} now implies that the maximal
submodule of the Verma module $M_{F}(\lambda ^{2})$ is generated by
five singular vectors with weights
\[
r_{\delta - \epsilon_1 -\epsilon_{3}}.\lambda ^{2}, \
r_{\alpha_{2}}.\lambda ^{2}, \
r_{\alpha_{3}}.\lambda ^{2}, \ r_{\alpha_{4}}.\lambda ^{2}, \
r_{\epsilon_{2}}.\lambda ^{2}.
\]
One can easily check that $f_{\frac{1}{2}( \epsilon_{1} - \epsilon_{2}
- \epsilon_{3} - \epsilon_{4})}(0)v^{(2)}$ is a singular vector with
weight $r_{\alpha_{4}}.\lambda ^{2}$, that
$f_{\epsilon _{4}}(0)v^{(2)}$ is a singular vector with
weight $r_{\alpha_{3}}.\lambda ^{2}$, and that
$f_{\epsilon_{3} - \epsilon_{4}}(0)v^{(2)}$ is a singular vector
with weight $r_{\alpha_{2}}.\lambda ^{2}$. This implies that
\begin{eqnarray*}
&& f_{\frac{1}{2}( \epsilon_{1} - \epsilon_{2}
- \epsilon_{3} - \epsilon_{4})}(0)v^{(2)}=0, \\
&& f_{\epsilon _{4}}(0)v^{(2)}=0 \quad \mbox{and} \\
&& f_{\epsilon_{3} - \epsilon_{4}}(0)v^{(2)}=0
\end{eqnarray*}
holds in $L_{F}(-\frac{5}{2}, -\frac{3}{2} \omega _{1})$. The
claim of lemma now follows immediately. $\;\;\;\;\Box$

\begin{thm}
\begin{eqnarray*}
& (1) & L_{F}(-\frac{5}{2},-\frac{3}{2} \omega _{1})
\cong L_{B}(-\frac{5}{2}, -\frac{3}{2} \bar{\omega} _{2}) \oplus
L_{B}(-\frac{5}{2}, -\frac{5}{2} \bar{\omega} _{2}+ \bar{\omega} _{4}) \\
& (2) & L_{F}(-\frac{5}{2}, -\frac{1}{2} \omega _{1}
-\frac{1}{2} \omega _{2}) \cong L_{B}(-\frac{5}{2}, -\frac{1}{2} \bar{\omega} _{1}
- \frac{1}{2} \bar{\omega} _{2} - \frac{1}{2} \bar{\omega} _{3}) \\
&& \quad \oplus
L_{B}(-\frac{5}{2}, -\frac{3}{2} \bar{\omega} _{1}+ \frac{1}{2} \bar{\omega} _{2}
- \frac{3}{2} \bar{\omega} _{3}+ \bar{\omega} _{4})  \\
& (3) & L_{F}(-\frac{5}{2}, -\frac{3}{2} \omega _{2}
+ \omega _{3}) \cong L_{B}(-\frac{5}{2}, -\frac{3}{2} \bar{\omega} _{1}
- \frac{1}{2} \bar{\omega} _{3}) \\
&& \quad \oplus
L_{B}(-\frac{5}{2}, -\frac{1}{2} \bar{\omega} _{1}-
\frac{3}{2} \bar{\omega} _{3} + \bar{\omega} _{4})
\end{eqnarray*}
\end{thm}
{\bf Proof:} It follows from relation (\ref{rel.lowest.conf.w})
that the lowest conformal weight of
$L_{F}(-\frac{5}{2},0)$-modules $L_{F}(-\frac{5}{2},-\frac{3}{2} \omega
_{1})$, $L_{F}(-\frac{5}{2}, -\frac{1}{2} \omega _{1}
-\frac{1}{2} \omega _{2})$ and $L_{F}(-\frac{5}{2}, -\frac{3}{2} \omega _{2}
+ \omega _{3})$ is $-\frac{3}{2}$. Let's prove (1).
Let $v_{\lambda}$ be any singular vector for $\hat{\frak g}_{B}$ in
$L_{F}(-\frac{5}{2},-\frac{3}{2} \omega _{1})$, with weight $\lambda \in
\hat{\frak h}^{*}$. $L_{B}(\lambda)$ is then a subquotient of
weak $L_{B}(-\frac{5}{2},0)$-module
$L_{F}(-\frac{5}{2},-\frac{3}{2} \omega _{1})$, which implies that
$L_{B}(\lambda)$ is an irreducible weak $L_{B}(-\frac{5}{2},0)$-module.
Thus, $L_{B}(\lambda)$ is isomorphic to some weak module listed in
Proposition \ref{p-B-n=1}. Furthermore,
Theorem \ref{prop.jedn.Vir} and Lemma \ref{lem.lowest.conf.w}
imply that $L_{B}(\lambda)$ is isomorphic to some weak module from
Lemma \ref{lem.lowest.conf.w} (3), and that
$v_{\lambda}$ is in the lowest conformal weight subspace of
$L_{F}(-\frac{5}{2},-\frac{3}{2} \omega _{1})$.
Using Lemma \ref{lem.7.2} (1), one can easily check that
$v^{(2)}$ and $f_{\frac{1}{2}( \epsilon_{1} + \epsilon_{2}
- \epsilon_{3} - \epsilon_{4})}(0)v^{(2)}$ are singular vectors
for $\hat{\frak g}_{B}$ in $L_{F}(-\frac{5}{2},-\frac{3}{2} \omega _{1})$
with weights $-\frac{5}{2} \Lambda _{0} -\frac{3}{2} \bar{\omega} _{2}$
and $-\frac{5}{2} \Lambda _{0} -\frac{5}{2} \bar{\omega} _{2}+ \bar{\omega}
_{4}$, respectively. Moreover, one can easily check that there are no singular vectors for
$\hat{\frak g}_{B}$ in
$L_{F}(-\frac{5}{2},-\frac{3}{2} \omega _{1})$ with weights
$\bar{\lambda} _{1}=-\frac{5}{2} \Lambda _{0}
-\frac{1}{2} \bar{\omega} _{1}
- \frac{1}{2} \bar{\omega} _{2} - \frac{1}{2} \bar{\omega} _{3}$,
$\bar{\lambda} _{2}=-\frac{5}{2} \Lambda _{0}-\frac{3}{2} \bar{\omega} _{1}+ \frac{1}{2} \bar{\omega} _{2}
- \frac{3}{2} \bar{\omega} _{3}+ \bar{\omega} _{4}$,
$\bar{\lambda} _{3}=-\frac{5}{2} \Lambda _{0}-\frac{3}{2} \bar{\omega} _{1}
- \frac{1}{2} \bar{\omega} _{3}$ and
$\bar{\lambda} _{4}=-\frac{5}{2} \Lambda _{0}-\frac{1}{2} \bar{\omega} _{1}-
\frac{3}{2} \bar{\omega} _{3} + \bar{\omega} _{4}$, since
$-\frac{5}{2} \Lambda _{0}-\frac{3}{2} \omega _{1}- \bar{\lambda} _{i}$
is not a sum of positive roots for $\frak g _{F}$, for
$i=1,2,3,4$.

Thus, $v^{(2)}$ and $f_{\frac{1}{2}( \epsilon_{1} + \epsilon_{2}
- \epsilon_{3} - \epsilon_{4})}(0)v^{(2)}$ are only singular vectors
for $\hat{\frak g}_{B}$ in $L_{F}(-\frac{5}{2},-\frac{3}{2} \omega _{1})$
which implies that
\begin{eqnarray*}
&& U(\hat{\frak g}_{B}). v^{(2)} \cong L_{B}(-\frac{5}{2}, -\frac{3}{2} \bar{\omega} _{2})
 \quad \mbox{and} \\
&& U(\hat{\frak g}_{B}). f_{\frac{1}{2}( \epsilon_{1} + \epsilon_{2}
- \epsilon_{3} - \epsilon_{4})}(0)v^{(2)} \cong L_{B}(-\frac{5}{2},
-\frac{5}{2} \bar{\omega} _{2}+ \bar{\omega} _{4}),
\end{eqnarray*}
and that
\begin{eqnarray*}
& & L_{F}(-\frac{5}{2},-\frac{3}{2} \omega _{1})
\cong L_{B}(-\frac{5}{2}, -\frac{3}{2} \bar{\omega} _{2}) \oplus
L_{B}(-\frac{5}{2}, -\frac{5}{2} \bar{\omega} _{2}+ \bar{\omega}
_{4}).
\end{eqnarray*}
Thus, we have proved claim (1).
Using Lemmas \ref{lem.lowest.conf.w} and \ref{lem.7.2}, one can
similarly obtain that
\begin{eqnarray*}
&& U(\hat{\frak g}_{B}). v^{(3)} \cong L_{B}(-\frac{5}{2}, -\frac{1}{2} \bar{\omega} _{1}
- \frac{1}{2} \bar{\omega} _{2} - \frac{1}{2} \bar{\omega} _{3})
 \quad \mbox{and} \\
&& U(\hat{\frak g}_{B}). f_{\frac{1}{2}( \epsilon_{1} - \epsilon_{2}
+ \epsilon_{3} - \epsilon_{4})}(0)v^{(3)} \cong
L_{B}(-\frac{5}{2}, -\frac{3}{2} \bar{\omega} _{1}+ \frac{1}{2} \bar{\omega} _{2}
- \frac{3}{2} \bar{\omega} _{3}+ \bar{\omega} _{4}),
\end{eqnarray*}
and that
\begin{eqnarray*}
& & L_{F}(-\frac{5}{2}, -\frac{1}{2} \omega _{1}
-\frac{1}{2} \omega _{2}) \cong L_{B}(-\frac{5}{2}, -\frac{1}{2} \bar{\omega} _{1}
- \frac{1}{2} \bar{\omega} _{2} - \frac{1}{2} \bar{\omega} _{3}) \\
&& \quad \oplus
L_{B}(-\frac{5}{2}, -\frac{3}{2} \bar{\omega} _{1}+ \frac{1}{2} \bar{\omega} _{2}
- \frac{3}{2} \bar{\omega} _{3}+ \bar{\omega} _{4})
\end{eqnarray*}
which implies claim (2). Furthermore
\begin{eqnarray*}
&& U(\hat{\frak g}_{B}). v^{(4)} \cong L_{B}(-\frac{5}{2}, -\frac{1}{2} \bar{\omega} _{1}-
\frac{3}{2} \bar{\omega} _{3} + \bar{\omega} _{4}),  \\
&& U(\hat{\frak g}_{B}). f_{\frac{1}{2}( \epsilon_{1} - \epsilon_{2}
- \epsilon_{3} + \epsilon_{4})}(0)v^{(4)} \cong
L_{B}(-\frac{5}{2}, -\frac{3}{2} \bar{\omega} _{1}
- \frac{1}{2} \bar{\omega} _{3}),
\end{eqnarray*}
and
\begin{eqnarray*}
& & L_{F}(-\frac{5}{2}, -\frac{3}{2} \omega _{2}
+ \omega _{3}) \cong L_{B}(-\frac{5}{2}, -\frac{3}{2} \bar{\omega} _{1}
- \frac{1}{2} \bar{\omega} _{3}) \\
&& \quad \oplus
L_{B}(-\frac{5}{2}, -\frac{1}{2} \bar{\omega} _{1}-
\frac{3}{2} \bar{\omega} _{3} + \bar{\omega} _{4})
\end{eqnarray*}
which implies claim (3). $\;\;\;\;\Box$

\section{Modules for $L_{F}(n-\frac{7}{2},0)$, $n \in \N$}
\label{sect.visi.nivoi}

In this section we study the category of
$L_{F}(n-\frac{7}{2},0)$-modules, for $n \in \N$.
Using Corollary \ref{cor.ulaganje-pol1} and
a certain polynomial from \cite{P}, one can easily obtain the classification of
irreducible objects in that category. Moreover, one can show that
this category is semisimple. We omit the
proofs in this section because of their similarity
with the proofs given in \cite{P}.

The next lemma follows from \cite[Lemma 18]{P} and Corollary \ref{cor.ulaganje-pol1}:

\begin{lem}\label{l.2.4.2}
Let
\begin{eqnarray*}
&\!\!\!\!\!\!& \!\! q(h)= \!\!\!\!\!\!\! {\displaystyle
\sum_{{(k_{1},\ldots ,k_{l}) \in {\Z}_{+}^{l} \atop \sum k_{i}=n}}} \frac{1}{k_{1}!4^{k_{1}}} \cdot
(h_{\epsilon_1}-2k_{2}-\ldots -2k_{l})\cdot \ldots \cdot
(h_{\epsilon_1}-2n+1)\cdot \\
& & \ \ \ \ \ \ \ \ \ \cdot (h_{\epsilon_1 - \epsilon_l}
-k_{l-1}-\ldots -k_{2}) \cdot \ldots \cdot (h_{\epsilon_1 - \epsilon_l}
-k_{l-1}-\ldots -k_{2}-k_{l}+1) \cdot \\
& & \ \ \ \ \ \ \ \ \ \cdot \ldots \cdot
h_{\epsilon_1 - \epsilon_2} \cdot \ldots \cdot (h_{\epsilon_1 -
\epsilon_2}-k_{2}+1).
\end{eqnarray*}
Then $q \in {\mathcal P}_{0}^{F}$.
\end{lem}

It can easily be checked (as in \cite[Lemma 21]{P}) that relation $q(\mu)=0$ for
$\mu \in P_{+}^{F}$ implies $( \mu , \epsilon_{1} ) \leq
n-\frac{1}{2}$. Using results from Subsection \ref{subsec.3.3} we get

\begin{lem}
Assume that $V_{F}(\mu), \mu \in P_{+}^{F}$ is an
$A(L_{F}(n-\frac{7}{2},0))$-module. Then
$( \mu , \epsilon_{1} ) \leq n-\frac{1}{2}$.
\end{lem}

Moreover, using weight arguments as in \cite[Lemma 22]{P} one can
prove the converse:

\begin{lem} Let $\mu \in P_{+}^{F}$, such that
$( \mu , \epsilon_{1} ) \leq n-\frac{1}{2}$. Then
$V_{F}(\mu)$ is an $A(L_{F}(n-\frac{7}{2},0))$-module.
\end{lem}

We get:

\begin{prop} The set
\begin{eqnarray*}
\left\{ V_{F}(\mu) \ \vert \ \mu \in P_{+}^{F}, \ ( \mu , \epsilon_{1} ) \leq
n-\frac{1}{2} \right\}
\end{eqnarray*}
provides the complete list of irreducible finite-dimensional
$A(L_{F}(n-\frac{7}{2},0))$-modules.
\end{prop}

It follows from Zhu's theory that:

\begin{thm}
The set
\begin{eqnarray*}
\left\{ L_{F}(n-\frac{7}{2}, \mu) \ \vert \ \mu \in P_{+}^{F},
\ ( \mu , \epsilon_{1} ) \leq
n-\frac{1}{2} \right\}
\end{eqnarray*}
provides the complete list of irreducible
$L_{F}(n-\frac{7}{2},0)$-modules.
\end{thm}

Furthermore, one can easily check (as in Lemma \ref{l.3.4}) that highest
weights $\lambda$ of these modules are admissible such that
\[
\hat{\Pi}^{\vee}_{\lambda}= \{ (\delta - \epsilon_{1})^{\vee},
\alpha_{1}^{\vee},\alpha_{2}^{\vee}, \alpha_{3}^{\vee} ,
\alpha_{4}^{\vee} \}.
\]

Using Proposition \ref{t.KW2} (as in \cite[Lemma 26, Theorem 27]{P}), one can easily
obtain:

\begin{thm}
Let $M$ be a $L_{F}(n-\frac{7}{2},0)$-module. Then
$M$ is completely reducible.
\end{thm}

\bibliography{thesis}

\begin{thebibliography}{FLM1}

\bibitem [A1]{A1}
D. Adamovi\'{c}, Some rational vertex algebras, { \it Glas. Mat. Ser.
III} {\bf 29}({\bf 49}) (1994), 25--40.

\bibitem [A2]{A2}
D. Adamovi\'{c}, Representations of vertex algebras associated
to symplectic affine Lie algebra at half-integer levels (in
Croatian), Ph.D. Thesis, University of Zagreb, 1996.

\bibitem [A3]{A3}
D. Adamovi\'{c}, A construction of admissible $A\sp {(1)}\sb 1$-modules of
level $-\frac 43$, { \it J. Pure Appl. Algebra} {\bf 196} (2005), 119--134.

\bibitem [AGO]{AGO}
R. C. Arcuri, J. F. Gomes and D. I. Olive,
Conformal subalgebras and symmetric spaces,
{ \it Nuclear Phys. B} {\bf 285} (1987), 327--339.


\bibitem [AM]{AM}
D. Adamovi\'{c} and A. Milas, Vertex operator algebras associated to
modular invariant representations for $A_{1}^{(1)}$, { \it Math.
Res. Lett.} {\bf 2} (1995), 563--575.


\bibitem[Ba]{Ba}
J. C. Baez, The octonions, {\it Bull. Amer. Math. Soc. (N.S.)} {\bf 39} (2002),
145--205.

\bibitem[Bo]{Bo}
R. E. Borcherds, Vertex algebras, Kac-Moody algebras, and the
Monster, {\it Proc. Natl. Acad. Sci. USA} {\bf 83} (1986),
3068--3071.


\bibitem[BB]{BB}
F. A. Bais and P. G. Bouwknegt,
A classification of subgroup truncations of the bosonic string,
{\it Nuclear Phys. B} {\bf 279} (1987), 561--570.


\bibitem[DL]{DL}
C. Dong and J. Lepowsky, {\em Generalized Vertex Algebras and Relative
Vertex Operators}, Progress in Math., {\bf Vol. 112}, Birkhauser,
Boston, 1993.

\bibitem[DLM1]{DLM1}
C. Dong, H. Li and G. Mason, Vertex operator algebras associated
to admissible representations of $\widehat{\rm sl}_ 2$,
{\it Comm. Math. Phys.} {\bf 184} (1997), 65--93.

\bibitem[DLM2]{DLM2}
C. Dong, H. Li and G. Mason, Simple currents and extensions of
vertex operator algebras, {\it Comm. Math. Phys.} {\bf 180} (1996), 671--707.


\bibitem[FHL]{FHL}
I. Frenkel, Y.-Z. Huang and J. Lepowsky, On axiomatic approaches to
vertex operator algebras and modules, {\em Mem. Amer. Math. Soc.}
{\bf 104}, 1993.


\bibitem[FLM]{FLM}
I. Frenkel, J. Lepowsky and A. Meurman, {\it Vertex Operator
Algebras and the Monster}, Pure and Appl. Math., {\bf Vol. 134},
Academic Press, Boston, 1988.

\bibitem[FM]{FM}
B. Feigin and F. Malikov, Fusion algebra at a rational level and cohomology of
nilpotent subalgebras of $\widehat{\rm sl}(2)$,
{\it Lett. Math. Phys.} {\bf 31} (1994), 315--325

\bibitem[FZ]{FZ}
I. Frenkel and Y.-C. Zhu, Vertex operator algebras associated to
representations of affine and Virasoro algebras, {\it Duke Math.
J.} {\bf 66} (1992), 123--168.

\bibitem[GKRS]{GKRS}
B. Gross, B. Kostant, P. Ramond and S. Sternberg, The Weyl character
formula, the half-spin representations, and equal rank subgroups,
{\it Proc. Natl. Acad. Sci. USA} {\bf 95} (1998), 8441--8442.


\bibitem[GPW]{GPW}
A. Ch. Ganchev, V. B. Petkova and G. M. T. Watts, A note on decoupling conditions
for generic level $\widehat{\rm sl}(3)\sb k$ and fusion rules,
{\em Nuclear Phys.} {\bf B571} (2000), 457--478.


\bibitem[K]{K}
V. G. Kac, {\em Infinite dimensional Lie algebras}, 3rd ed.,
Cambridge Univ. Press, Cambridge, 1990.


\bibitem[KW1]{KW1}
V. Kac and M. Wakimoto, Modular invariant representations of infinite dimensional Lie algebras
and superalgebras, {\em Proc. Natl. Acad. Sci. USA} {\bf 85} (1988), 4956--4960.

\bibitem[KW2]{KW2}
V. Kac and M. Wakimoto, Classification of modular invariant representations
of affine algebras, in {\em Infinite Dimensional Lie algebras and groups},
Advanced Series in Math. Phys. {\bf 7}, World Scientific,
Teaneck NJ, 1989.

\bibitem[L]{L}
H.-S. Li, Local systems of vertex operators, vertex superalgebras
and modules. {\em J. Pure Appl. Algebra} {\bf 109} (1996), 143--195.


\bibitem[LL]{LL} J. Lepowsky and H. Li, {\it Introduction to vertex operator algebras
and their representations}, Progress in Math., {\bf Vol. 227}, Birkhauser,
Boston, 2004.


\bibitem[MP]{MP}
A. Meurman and M. Primc, {\em Annihilating fields of standard
modules of $sl(2, \C)\tilde{}$ and combinatorial identities}, Mem.
Amer. Math. Soc. {\bf 137}, AMS, Providence RI, 1999.

\bibitem[P]{P}
O. Per\v{s}e, Vertex operator algebras
associated to type $B$ affine Lie algebras
on admissible half-integer levels,
{\em J. Algebra} {\bf 307} (2007), 215--248.

\bibitem[PR]{PR}
T. Pengpan and P. Ramond, M(ysterious) paterns in ${\rm SO} (9)$.
Looking forward: frontiers in theoretical science (Los Alamos, NM,
1998), {\em Phys. Rep.} {\bf 315} (1999), 137--152.

\bibitem[SW]{SW}
A. N. Schellekens and N. P. Warner, Conformal subalgebras of
Kac-Moody algebras, {\em Phys. Rev. D} {\bf 34} (1986), 3092--3096.


\bibitem[W]{W}
M. Wakimoto, {\em Lectures on infinite-dimensional Lie algebra},
World Scientific, River Edge NJ, 2001.

\bibitem[Z]{Z}
Y.-C. Zhu, Modular invariance of characters of
vertex operator algebras, {\em J.
Amer. Math. Soc.} {\bf 9} (1996), 237--302.


\end{thebibliography}
\bibliographystyle{plain}

\vskip 1cm

Department of Mathematics, University of Zagreb, Bijeni\v{c}ka 30,
\linebreak 10000 Zagreb, Croatia

E-mail address: perse@math.hr

\end{document}